\documentclass{amsproc}

\usepackage{graphics}

\newtheorem{theorem}{Theorem}[section]
\newtheorem{lemma}[theorem]{Lemma}
\newtheorem{corollary}[theorem]{Corollary}
\newtheorem{proposition}[theorem]{Proposition}
\newtheorem{definition}[theorem]{Definition}

\newtheorem{remark}[theorem]{Remark}

\numberwithin{equation}{section}

\newcommand{\bz}{{\mathbb B}}
\newcommand{\cz}{{\mathbb C}}
\newcommand{\hz}{{\mathbb H}}
\newcommand{\nz}{{\mathbb N}}
\newcommand{\rz}{{\mathbb R}}

\newcommand{\calC}{\mathcal{C}}
\newcommand{\calD}{\mathcal{D}}
\newcommand{\calH}{\mathcal{H}}
\newcommand{\calK}{\mathcal{K}}
\newcommand{\calL}{\mathcal{L}}
\newcommand{\calM}{\mathcal{M}}
\newcommand{\calO}{\mathcal{O}}
\newcommand{\calS}{\mathcal{S}}

\newcommand{\amax}{A_{\text{\rm max}}}
\newcommand{\amin}{A_{\text{\rm min}}}
\newcommand{\ci}{\mathcal{C}^\infty}
\newcommand{\cicomp}{\mathcal{C}^\infty_{\text{\rm comp}}}
\newcommand{\cigb}{\mathcal{C}^{\infty,\gamma}({\mathbb B})}
\newcommand{\dbar}{d\hspace*{-0.08em}\bar{}\hspace*{0.1em}}
\newcommand{\dvz}{\,\mbox{\rm div}\,} 
\newcommand{\eps}{\varepsilon}
\newcommand{\grad}{\,\mbox{\rm grad}\,}

\newcommand{\hsgpb}{\mathcal{H}^{s,\gamma}_p({\mathbb B})}

\newcommand{\im}{\text{\rm Im}\,}
\newcommand{\intb}{\text{\rm int}\,{\mathbb B}}
\newcommand{\op}{\text{\rm op}}
\newcommand{\opm}[1]{\text{\rm op}_M^{#1}}
\newcommand{\pit}{\,{\widehat{\otimes}}_\pi\,}
\newcommand{\re}{\text{\rm Re}\,}
\newcommand{\rpbar}{\overline{{\mathbb R}}_+}
\newcommand{\skp}[2]{\langle#1,#2\rangle}
\newcommand{\spk}[1]{\left<#1\right>}
\newcommand{\st}{\mbox{\boldmath$\;|\;$\unboldmath}}
\newcommand{\trinorm}[1]%
    {|\hspace*{-1pt}|\hspace*{-1pt}|#1|\hspace*{-1pt}|\hspace*{-1pt}|}

\renewcommand{\Re}{{\rm Re}\,}
\renewcommand{\Im}{{\rm Im}\,}
        
\begin{document}

\title[Bounded Imaginary Powers of Cone Differential Operators]
      {Bounded Imaginary Powers of Differential Operators\\ 
      on Manifolds with Conical Singularities}
\author{S.\ Coriasco$^\dagger$}
\address{Universit\'a di Torino, Dipartimento di Matematica,
         Via Carlo Alberto 10, 10123 Torino, Italy}
\email{coriasco@dm.unito.it}
\author{E.\ Schrohe}
\address{Universit\"at Potsdam, Institut f\"ur Mathematik, 
         Postfach 60 15 53, 14415 Potsdam, Germany}
\email{schrohe@math.uni-potsdam.de}
\author{J.\ Seiler}
\address{Universit\"at Potsdam, Institut f\"ur Mathematik, 
         Postfach 60 15 53, 14415 Potsdam, Germany}
\email{seiler@math.uni-potsdam.de}
\date{\today}
\keywords{Complex powers, Manifolds with conical singularities, 
Resolvent.\newline
\mbox{\hspace{0.65cm}$^\dagger$}Supported by the E.U. 
          Research and Training Network ``Geometric Analysis''.}
\begin{abstract}
 We study the minimal and maximal closed extension of a differential 
 operator $A$ on a manifold $B$ with conical singularities, when $A$ 
 acts as an 
 unbounded operator on weighted $L_p$-spaces over $B$, $1<p<\infty$. 
 
 Under suitable ellipticity assumptions we can define a family of complex 
 powers $A^z$, $z\in\cz$. We also obtain sufficient information on the resolvent 
 of $A$ to show the boundedness of the purely imaginary powers. 
 
 Examples concern unique solvability and maximal regularity for the solution 
 of the Cauchy problem for the Laplacian on 
 conical manifolds as well as certain quasilinear diffusion equations. 
\end{abstract}

\maketitle

\markboth{\uppercase{Bounded Imaginary Powers of Cone Differential Operators}}
         {\uppercase{S.\ Coriasco, E.\ Schrohe, J.\ Seiler}}

\tableofcontents

\setlength{\parskip}{5pt}

\section{Introduction}
Seeley's classical paper \cite{Seel0}, published in 1967, showed in a 
striking way how pseudodifferential techniques could be applied to analyze 
complex powers of elliptic (pseudo-)differential operators on closed 
manifolds. Replacing the resolvent in the Dunford integral by a 
parameter-dependent parametrix, he obtained a representation of the powers
that was precise enough to deduce a wealth of information on eigenvalue 
asymptotics, zeta functions, and index theory. Seeley also extended his 
results to differential boundary value problems. In 1971 he showed the 
boundedness of the purely imaginary powers on $L_p$-spaces, \cite{Seel2}. 

At that time the principal motivation for these studies was the 
description of interpolation spaces. Additional interest in 
the behavior of imaginary powers came from Dore and Venni's 1987 article 
\cite{DoVe}, in which they showed how the boundedness of imaginary 
powers can be used to derive results on maximal regularity for evolution
equations.

Meanwhile, bounded imaginary powers or even the existence of a bounded 
$H^\infty$ calculus \cite{Mci} have been established in many situations,
e.g.\ in abstract settings \cite{PrSo2}, 
for classes of differential operators on $\rz^n$ and smooth manifolds 
\cite{AHS}, boundary value problems on bounded and certain unbounded 
domains in $\rz^n$, \cite{Duon}, \cite{PrSo}, \cite{SoTh}, as well as for
operators in Boutet de Monvel's calculus \cite{Sohr}. 

We shall focus here on the case of a manifold with conical singularities. 
This is a Hausdorff space,~$B$, that is a smooth manifold outside a 
finite number of singular points, while, close to each of these, it has the
structure of a cone with smooth, closed cross-section. Blowing up $B$ near 
its singular points, we obtain a manifold $\bz$ with boundary 
$\partial\bz=:X$. 

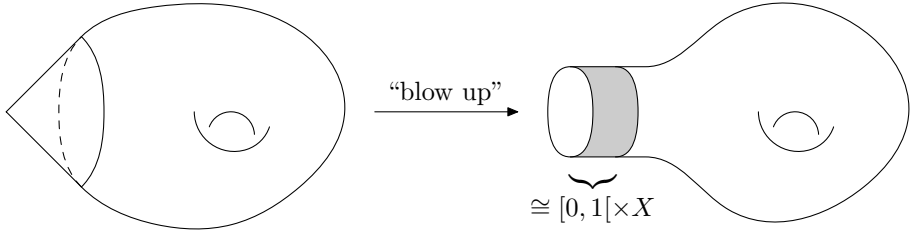
\begin{center}
    \includegraphics{figures.2}
\end{center}

Near the boundary, we fix a splitting of coordinates $(t,x)\in[0,1[\times 
X$. Rather than on $B$, the analysis will be performed on $\bz$ 
(respectively the interior of $\bz$). We consider so-called {\em cone} or 
{\em Fuchs-type} differential operators, i.e., operators which close to the 
boundary are of the form 
 \begin{equation}\label{intro1}
   A=t^{-\mu}\mathop{\mbox{\Large$\sum$}}_{j=0}^\mu 
   a_j(t)(-t\partial_t)^j,
 \end{equation}
where each $a_j\in\ci(\rpbar,\text{\rm Diff}^{\mu-j}(X))$ is a smooth 
family of differential operators on the cross-section. Such an $A$ acts as 
an unbounded operator 
$A:\cicomp(\intb)\subset\calH^{0,\gamma}_p(\bz)\to\calH^{0,\gamma}_p(\bz)$, 
where the space $\calH^{0,\gamma}_p(\bz)$ away from the boundary coincides 
with $L_p(\bz)$ and near the boundary with 
 $$t^{\gamma-\frac{n+1}{2}}L_p([0,1[\times X,\mbox{$\frac{dt}{t}$}dx),\qquad 
   n=\dim X.$$
Here, $1<p<\infty$, and $\gamma$ is an arbitrary real number. 
Justified by the fact that a change to polar coordinates shows the
equivalence
$L_p(\rz^{n+1})=t^{-\frac{n+1}{p}}L_p(\rz_+\times
S^n,\mbox{$\frac{dt}{t}$}d\varphi)$ 
we define the space $L_p(B)$ as $\calH_p^{0,\gamma_p}(\bz)$ for 
$\gamma_p=(n+1)(\frac{1}{2}-\frac{1}{p})$.

Let $\Lambda_\Delta=\Lambda_\Delta(\theta)$ denote a closed sector in the 
complex plane, symmetric about the negative real half-axis and of 
aperture $2(\pi-\theta)$ for some $0<\theta<\pi$. We find conditions  
(Definition \ref{ellipticity}) on $A$, which depend on $\gamma\in\rz$ but 
not on $1<p<\infty$, that ensure the following: 
 \begin{itemize}
  \item[i)] the {\em closure} $\amin$ of $A$ has no spectrum in 
   $\Lambda_\Delta\cap\{|\lambda|>R\}$;
  \item[ii)] the resolvent satisfies the uniform estimate 
   $\|(\lambda-A)^{-1}\|_{\calL(\calH^{0,\gamma}_p(\bz))}\le 
   c_p\,|\lambda|^{-1}$. 
 \end{itemize} 
Moreover, we obtain very precise information on the structure of the 
resolvent. For this and i) see Theorem \ref{parametrix}; ii) is shown in 
Proposition \ref{estimate}. We also give conditions (Remark 
\ref{ellipticity2}) implying that the {\em maximal extension} $\amax$ 
satisfies statements analogous to i) and ii). The symmetry 
about the negative real axis, which here plays the role of a
ray of minimal growth in the sense of Agmon \cite{Agmo}, 
\cite{Seel2}, is not essential.
The case of an arbitrary symmetry axis $\{te^{i\varphi}\st t\ge0\}$ can be
reduced to our situation, replacing $A$ by $e^{-i\varphi}A$. 

Since $\amin$ (respectively $\amax$) in the above case has compact 
resolvent, any ``keyhole'' region $\Lambda$, consisting 
of the sector $\Lambda_\Delta$ and an arbitrarily small ball around 
zero, only contains finitely many elements 
of the spectrum. Assuming that zero is the only spectral point in the 
keyhole (or, alternatively, shrinking the angle of the sector and possibly 
rotating $A$ a little), we define {\em complex powers} $A^z_{\min}$ 
(respectively $A^z_{\max}$) for all complex $z$ with negative real part. 
This is done in terms of a Dunford integral, integrating the resolvent 
against $\lambda^z$ along the boundary of the keyhole. Using the specific 
structure of the resolvent, we show (Theorem \ref{bip}) that 
 \begin{itemize}
  \item[iii)] $\|A^z_{\min}\|_{\calL(\calH^{0,\gamma}_p(\bz))}\le 
   c_p\,e^{|\im z|\theta}$ uniformly for all $z$ with $|\re z|$ 
   sufficiently small 
 \end{itemize} 
as well as the analogous estimate for $A^z_{\max}$ (Theorem 
\ref{bip2}). Consequently, the purely imaginary powers $A^{iy}_{\min}$ 
(respectively $A^{iy}_{\max}$), $y\in\rz$, exist as suitable limits and 
satisfy an estimate as in iii). 

It should be noted that both the construction of the complex powers 
and the boundedness of the imaginary powers only rely on the 
information about the resolvent provided by Theorem \ref{parametrix}. 
Our conclusions therefore carry over to all situations where the 
resolvent has this structure.

The key to the above described results is, similar to Seeley's classical 
concept, to view $\lambda-A$ as an element of a calculus of 
parameter-dependent pseudodifferential operators on $\bz$, and to express 
$(\lambda-A)^{-1}$ within this calculus. In our context, the appropriate 
calculus is Schulze's parameter-dependent cone algebra, cf.\  
for example \cite{Schu1}, \cite{EgSc}. 
 The conditions we impose on $A$ are,
more or less, ellipticity conditions on $\lambda-A$ within this calculus.
We require three associated objects not to have
spectrum in the sector $\Lambda_\Delta$. The first is the usual homogeneous
principal symbol of $A$, defined on the cotangent bundle over the interior
of $\bz$. The second is the so-called {\em rescaled symbol}, which reflects
the behavior of the principal symbol near the boundary. The third is the
so-called {\em model cone} operator $\widehat{A}$, which acts as an
unbounded operator in Sobolev spaces on the infinite cylinder $\rz_+\times
X$. It is induced by freezing the coefficients of $A$ at the boundary, i.e.,
using the notation from \eqref{intro1}, $\widehat{A}=t^{-\mu}\sum_{j=0}^\mu
a_j(0)(-t\partial_t)^j$. 

In order to separate the more general functional-analytic issues from 
the specific difficulties related to conical singularities, 
we give a review of several basic facts about complex powers of unbounded
operators on a Banach space in Section \ref{semigroup}, while in Section 
\ref{diffop} we briefly discuss Fuchs-type operators. Sections \ref{resolvent}, 
\ref{boimpo} and \ref{closed} are devoted to the proof of the results
stated above. 

In Section \ref{laplace} we treat an example and show how our work can be 
combined with that of Dore and Venni to obtain results on unique solvability
and maximal regularity for the non-homogeneous Cauchy problem in $L_p(B)$:
 $$\dot u(\tau)-\Delta u(\tau)=f(\tau),\qquad u(0)=0.$$ 
Here, $\Delta$ is the Laplace-Beltrami operator for a Riemannian metric with 
a conical degeneracy, and we assume $\dim \bz>4$. 
Based on this observation, 
we consider  in Section \ref{sec:quasi} quasilinear diffusion equations of
the form 
$$\dot u(\tau) - \dvz(a(t^cu)\grad u)(\tau) = f(u,\tau)+g(\tau), \qquad u(0) =
u_0,$$
in weighted $L_p$-spaces on $\bz$, and show unique solvability  with the help
of  an abstract result by Cl\'ement and Li
\cite{ClLi}.
Here, $a$, the diffusion coefficient, is a smooth positive function; 
we assume that it depends on $t^cu$, where $t$ is is a smooth function on
$\bz$, coinciding with the distance to the boundary near $\partial \bz$, and
$c$ is a positive constant.

It is clear that physically relevant applications would require our
understanding of the Laplacian on lower-dimensional
manifolds and of boundary value problems. 
Both topics are presently under investigation. 
As their analysis, however, is considerably more
complicated, it makes sense to focus first on 
the present situation, where the ideas and techniques
can be explained more easily. 

An appendix relates the structure of the resolvent as 
we use it to that given in earlier work by Schulze and that of Gil 
\cite{Gil}. Moreover, we collect a few definitions and 
notions in Section \ref{notation}. 

\vspace{0.2cm}

{\sc Acknowledgement}: We thank M.\ Korey (Potsdam) and M.\ Hieber
(Darmstadt) for several valuable discussions.

\section{Complex powers of operators in a Banach space}\label{semigroup}
Let us recall some well-known facts on complex powers of a closed, densely 
defined operator 
 \begin{equation*}
  A:\calD(A)\subset F\longrightarrow F
 \end{equation*}
in a Banach space $F$, cf.\ for example \cite{Seel2}. We denote by 
$\Lambda=\Lambda(\delta,\theta)$ the keyhole region
 $$\Lambda(\delta,\theta)=\{\lambda\in\cz\st|\lambda|\le\delta
                           \text{ or }|\arg\lambda|\ge\theta\}$$
with $\delta>0$ and $0<\theta<\pi$. We assume that 
 \begin{itemize}
  \item[(A1)] The spectrum of $A$ has empty intersection with 
   $\Lambda\setminus\{0\}$. 
  \item[(A2)] $\|(\lambda-A)^{-1}\|_{\calL(F)} |\lambda|$ is uniformly 
   bounded for large $\lambda\in\Lambda$. 
 \end{itemize} 
If $\calC$ is the parametrization of the boundary of $\Lambda$, cf.\ 
\eqref{not3}, we let 
 \begin{equation}\label{atoz}
  A^z=\frac{1}{2\pi i}\int_\calC\lambda^z(\lambda-A)^{-1}\,d\lambda,\qquad
  z\in\hz=\{z\in\cz\st\re z<0\}.
 \end{equation}
Here, $\lambda^z$ is defined via the logarithm 
$\log\lambda=\log|\lambda|+i\arg\lambda$, where $-\pi<\arg\lambda<\pi$. 
Since the integrand is $O(|\lambda|^{-1+\re z})$, the integral is 
absolutely convergent, and thus \eqref{atoz} defines continuous operators 
$A^z\in\mathcal{L}(F)$. The notation $A^z$ should be viewed with a little 
care, since $A^{-1}$ in general is not the inverse of $A$, which is not 
required to exist. 
\begin{remark}\label{group}
 Under conditions \text{\rm (A1), (A2)}, the function 
 $z\mapsto A^z:\hz\to\calL(F)$ is holomorphic and satisfies the semi-group 
 property 
  $$A^zA^w=A^{z+w},\qquad z,w\in\hz.$$ 
\end{remark}
If one furthermore imposes that for some positive constant $c$ 
 \begin{itemize}
  \item[(A3)] $\|A^z\|_{\calL(F)}$ is uniformly bounded in the rectangle 
   $-c\le\re z<0$, $|\im z|\le k$ for any $k\in\nz$, 
 \end{itemize}
then the limits of $A^z$ for $z\to iy$ exist for any real $y$. More 
precisely: 
\begin{remark}\label{limit}
 Under conditions \text{\rm (A1), (A2)}, and \text{\rm (A3)}, the limits  
  $$A^{iy}f=\lim_{\hz\,\ni z\to iy}A^zf$$
 exist for any real number $y$ and any $f\in F$, and thus define 
 operators $A^{iy}\in\mathcal{L}(F)$. Furthermore $A^{iy}f=A^{-1+iy}Af$ for 
 $f\in\calD(A)$. In particular, if we set  
  $$E_0=\frac{1}{2\pi i}\int_{|\lambda|=\delta}
    (\lambda-A)^{-1}\,d\lambda,$$
 then $E_0$ is a projection in $F$ and $A^0=1-E_0$. 
\end{remark}
Remark \ref{limit} could be rephrased as follows: Under conditions 
\text{\rm (A1), (A2)}, and \text{\rm (A3)}, the operators 
 $$T^z:=A^z+E_0,\qquad z\in\hz,$$ 
form an analytic semi-group (with $\lim_{z\to 0}T^zf=f$ for any $f\in F$) 
and there exist constants $c\ge1$ and $\omega\ge0$ such that 
 $$\|T^z\|_{\mathcal{L}(F)}\le c\,e^{\omega|z|}, \qquad z\in\hz.$$ 
Moreover, $(1-E_0)A+E_0:\mathcal{D}(A)\to F$ is an isomorphism, whose 
inverse is $T^{-1}$. 

In concrete situations the problem is to analyze whether an operator $A$ 
satisfies conditions (A1), (A2), (A3), and then to find the best possible 
constant $\omega$. Fundamental works on this topic are due to Seeley 
\cite{Seel0}, \cite{Seel1}, \cite{Seel2}, where he gives criteria ensuring 
that a differential operator on a compact manifold (with boundary) has 
these properties. The main object of the present paper is to give such 
criteria for differential operators on manifolds with conical 
singularities. 


\section{Cone differential operators}\label{diffop} 
We consider a differential operator $A: \calC^\infty(\bz,E)\to
\calC^\infty(\bz,E)$ acting on sections of a vector bundle $E$ over $\bz$. 
We may assume that $E$ respects the product structure near the boundary
$\partial \bz = X$,
i.e.\ is the  pull-back of a vector bundle $E_X$ over $X$ under the
canonical projection $[0,1[\times X\to X$. The operator $A$ is called a 
{\em cone differential} or {\em Fuchs-type} operator, cf.\ \cite{EgSc}, 
\cite{Lesc}, if, near the boundary it is of the form 
 \begin{equation}\label{diffop1}
  A=t^{-\mu}\mathop{\mbox{\Large$\sum$}}_{j=0}^\mu a_j(t)(-t\partial_t)^j,
 \end{equation}
where $a_j\in\ci(\rpbar,\text{\rm Diff}^{\mu-j}(X; E_X,E_X))$
are functions, smooth 
up to the boundary, with values in the differential operators on $X$. 
In order to keep the notation simple, we shall not indicate the bundles and
write $\calC^\infty  (\bz)$, ${\rm Diff}^{\mu-j}(X)$, etc. 

We can 
rewrite \eqref{diffop1} as 
 \begin{equation}\label{diffop2}
  A=t^{-\mu}\opm{\gamma+\mu-\frac{n}{2}}(f),\qquad
  f(t,z)=\mathop{\mbox{\Large$\sum$}}_{j=0}^\mu a_j(t)z^j,
 \end{equation}
where the Mellin pseudodifferential operator is defined by 
 \begin{equation}\label{mellop}
  [\opm{\gamma+\mu-\frac{n}{2}}(f)u](t)=
  \int_{\re z=\frac{n+1}{2}-\gamma-\mu}t^{-z}f(t,z)(\calM u)(z)\,\dbar z, 
  \qquad u\in\cicomp([0,1[\times X).
 \end{equation}
Here, $\gamma\in\rz$ is arbitrary, and \eqref{mellop} is 
independent of the choice of $\gamma$. We keep $\gamma$ in the 
notation, since we shall consider extensions of $A$ to different weighted 
Sobolev spaces, the weight being given by $\gamma$:  
\begin{definition}\label{sobolev}
 For $s\in\nz_0$, $\gamma\in\rz$, and $1<p<\infty$ we introduce $\hsgpb$ as 
 the space of all functions $u\in H^s_{p,\text{\rm loc}}(\intb)$ such that 
  $$t^{\frac{n+1}{2}-\gamma}(t\partial_t)^k\partial_x^\alpha
    (\omega u)(t,x)\;\in\;
    L_p(\rz_+\times X,\mbox{$\frac{dt}{t}dx$})\qquad
    \forall\;k+|\alpha|\le s$$
 for some cut-off function $\omega\in\cicomp([0,1[)$. 
\end{definition}
Recall that a cut-off function is a function $\omega\in \cicomp(\intb)$ such that $\omega \equiv 1$ near $t=0$.
The definition of the Banach spaces $\hsgpb$ naturally extends to real $s$ as
well as to spaces of vector bundles over $\bz$.  For more details see Section
\ref{notation}. 

For any $s$, $\gamma$, and $p$, the operator $A$ induces continuous
mappings 
 $$
  A:\calH^{s+\mu,\gamma+\mu}_p(\bz)\longrightarrow\hsgpb. 
 $$
With $A$ we associate three symbols. The first is the usual 
{\em homogeneous principal symbol} 
$\sigma^{\mu}_\psi(A)\in\ci(T^*(\intb)\setminus0)$, taking values in the
corresponding bundle homomorphisms. In local  coordinates 
near the boundary
 \begin{equation}\label{principal}
  \sigma^{\mu}_\psi(A)(t,x,\tau,\xi)=t^{-\mu}
  \mathop{\mbox{\Large$\sum$}}_{j=0}^\mu 
  \sigma^{\mu-j}_\psi(a_j)(t,x,\xi)(-it\tau)^j
 \end{equation}
Dropping the factor $t^{-\mu}$, replacing $t\tau$ by $\tau$, and inserting 
$t=0$, we obtain 
 $$\widetilde{\sigma}^{\mu}_\psi(A)(x,\tau,\xi):=
   \mathop{\mbox{\Large$\sum$}}_{j=0}^\mu 
   \sigma^{\mu-j}_\psi(a_j)(0,x,\xi)(-i\tau)^j,$$ 
which yields the {\em rescaled symbol} of $A$, 
 \begin{equation}\label{principal2}
  \widetilde{\sigma}^{\mu}_\psi(A)\;\in\;\ci((T^*X\times\rz)\setminus0). 
 \end{equation}
The third one is the {\em conormal symbol}
 \begin{equation}\label{conormal}
  \sigma^\mu_M(A)(z)=f(0,z)=
   \mathop{\mbox{\Large$\sum$}}_{j=0}^\mu a_j(0)z^j,
 \end{equation}
a function of $z\in\cz$ with values in the differential operators on $X$. 
\begin{remark}\label{coneell}
 The operator $A$ is called elliptic with respect to the weight 
 $\gamma+\mu$, if both the homogeneous principal symbol and 
 the rescaled symbol are invertible on their respective domains, and 
  \begin{equation}\label{conormal2}
   \sigma^\mu_M(A)(z):H^s(X)\to H^{s-\mu}(X),\qquad
    \re z=\mbox{$\frac{n+1}{2}-\gamma-\mu$}, 
  \end{equation}
 is an isomorphism for all $z$ on that line.
 
 It can be shown, \cite{ScSe1}, 
 Theorem \text{\rm 3.13}, that $A$ is elliptic if and only if the operators 
 $A:\calH^{s+\mu,\gamma+\mu}_p(\bz)$ $\longrightarrow\hsgpb$
 are Fredholm for any $s$ and $p$. 
\end{remark}
We shall consider $A$ as the operator 
 \begin{equation}\label{extension}
  A:\calD(A)=\calH^{\mu,\gamma+\mu}_p(\bz)\subset
  \calH^{0,\gamma}_p(\bz)\longrightarrow\calH^{0,\gamma}_p(\bz). 
 \end{equation}
\begin{remark}
  \begin{itemize}
      \item[a)] In case $A$ is elliptic with respect to $\gamma+\mu$,
                \eqref{extension} is 
                the closure of $A$ considered on the domain $\cicomp(\intb)$. 
      \item[b)] By the spectral invariance of the cone algebra, 
                \cite{ScSe1}, Theorem \text{\rm 3.14}, the spectrum of $A$ is 
                independent of $1<p<\infty$. 
  \end{itemize}
\end{remark}
With $A$ we associate the {model cone operator}, which acts in 
Sobolev spaces on the infinite cylinder 
 \begin{equation}\label{xwedge}
  X^\wedge:=\rz_+\times X. 
 \end{equation}
 Let $(t,x)$ denote cylindrical coordinates on $X^\wedge$. Then
 $H^s_{p,\text{\rm cone}}(X^\wedge)$ is the space of all 
 distributions $u$ whose push-forward under conical coordinates 
 $(t,tx)$ belongs to $H^s_{p}(\rz^{1+n})$ $($for details see 
   \cite{ScSc}, Section {\rm 4.2}$)$.
\begin{definition}\label{ksgp}
 For $s,\gamma\in\rz$ and $1<p<\infty$ the spaces 
 $\calK^{s,\gamma}_p(X^\wedge)$ consist of all distributions $u\in 
 H^s_{p,\text{\rm loc}}(X^\wedge)$ satisfying, for some cut-off function 
 $\omega\in\cicomp([0,1[)$, 
  $$\omega u\in\hsgpb,\qquad
    (1-\omega)u\in H^s_{p,\text{\rm cone}}(X^\wedge).$$
\end{definition}
Freezing the coefficients of $A$ at $t=0$, we obtain the 
model cone operator $\widehat{A}$,  
 \begin{equation}\label{ahat}
  \widehat{A}=t^{-\mu}\mathop{\mbox{\Large$\sum$}}_{j=0}^\mu
                             a_j(0)(-t\partial_t)^j:
  \calK^{\mu,\gamma+\mu}_p(X^\wedge)\longrightarrow
  \calK^{0,\gamma}_p(X^\wedge).
 \end{equation}
\begin{remark}
 \begin{itemize}
  \item[a)] If $A$ is elliptic with respect to the weight $\gamma+\mu$ and 
   satisfies the ellipticity condition \text{\rm (E1)} introduced below, 
   then it can be shown that \eqref{ahat} is the closure of $\widehat{A}$ 
   considered on the domain $\cicomp(X^\wedge)$.  
  \item[b)] If we set $\widehat{a}(\lambda)=\lambda-\widehat{A}$ with 
   $\widehat{A}$ from \eqref{ahat}, then $\widehat{a}(\lambda)$ corresponds to 
   the so-called {\em principal edge symbol} of $\lambda-A$, if we view 
   $\lambda-A$ as a constant coefficient edge symbol in the framework of 
   Schulze's theory of pseudodifferential operators on manifolds with 
   edges, cf.\ for example \cite{EgSc}. 
 \end{itemize}
\end{remark} 
It is worth mentioning that $\widehat{a}(\lambda)=\lambda-\widehat{A}$ is a 
homogeneous function in a specific way. Namely if we define for $\varrho>0$
 \begin{equation}\label{kappa}
  \kappa_\varrho:\cicomp(X^\wedge)\to\cicomp(X^\wedge),\quad 
  (\kappa_\varrho u)(t,x)=\varrho^{\frac{n+1}{2}}u(\varrho t,x),
 \end{equation}
then these operators extend by continuity to isomorphisms in 
$\calL(\calK^{s,\gamma}_p(X^\wedge))$ and 
 \begin{equation}\label{twisted}
  \widehat{a}(\varrho^\mu\lambda)=  
  \varrho^\mu\,\kappa_\varrho\,\widehat{a}(\lambda)\,\kappa_\varrho^{-1}. 
 \end{equation}
In particular, $\text{\rm spec}(\widehat{A})$ is a closed conic subset of 
the complex plane. 

Near the  boundary of $\bz$ we can write 
 $\lambda-A=t^{-\mu}\opm{\gamma+\mu-\frac{n}{2}}(h)(\lambda)$ 
with the parameter-dependent Mellin symbol 
 \begin{equation}\label{tildeh}
  h(t,z,\lambda)=\tilde{h}(t,z,t^\mu\lambda),\qquad
  \tilde{h}(t,z,\lambda)=\lambda-f(t,z),
 \end{equation}
and $f$ from \eqref{diffop2}. 

\section{The resolvent of cone differential operators}\label{resolvent}
To describe the structure of the resolvent we 
recall some elements from the theory of parameter-dependent cone 
pseudodifferential operators, starting with
the smoothing remainders of the calculus. 
To this end we introduce a family of Fr\'echet spaces of smooth functions on 
$\intb$ and $X^\wedge$, respectively.  
\begin{definition}\label{test}
 For $\gamma\in\rz$ we let $\cigb$ denote the space of all $u \in 
 \ci(\intb)$ such that
  \begin{equation}\label{cig}
      \sup_{0<t<1}t^{\frac{n+1}{2}-\gamma}
    \trinorm{\log^l t~(t\partial_t)^k( u)(t,\cdot)}<\infty\qquad
    \forall\;k,l\in\nz_0
  \end{equation}
 for any semi-norm $\trinorm{\cdot}$ of $\ci(X)$. Similarly,
 $\calS^\gamma_0(X^\wedge)$ is the space of all $u \in 
 \ci(X^\wedge)$ which are rapidly decreasing as $t \to \infty$ and 
 satisfy \eqref{cig}.
\end{definition}
We shall say that an operator $G$ has a kernel $k$ with respect to the 
$\calH^{0,0}_2(\bz)$-scalar product if 
 $$(Gu)(y)=\skp{k(y,\cdot)}{\overline{u}}_{\calH^{0,0}_2(\bz)}=
   \int_\bz k(y,y')u(y')\,t(y')^ndy',\qquad u\in\cicomp(\intb),$$
where $t$ denotes a boundary defining function on $\bz$ and $dy'$ refers 
to a density on $2\bz$, the double of $\bz$. We shall 
use the analogous notion for operators on $X^\wedge$, based on the 
scalar product of $\calK^{0,0}_2(X^\wedge)=L_2(X^\wedge,t^ndtdx)$. 
\begin{definition}\label{green1}
 An operator-family $G=G(\lambda)$, $\lambda \in \Lambda$, belongs to 
 $C^{-\infty}_G(\bz;\Lambda,\gamma)$, $\gamma\in\rz$, if there exists 
 an $\eps=\eps(G)>0$ such that $G(\lambda)$ has a kernel 
 $k(\lambda)=k(\lambda,\cdot,\cdot)$ with 
 respect to the $\calH^{0,0}_2(\bz)$-scalar product and 
  $$k(\lambda,y,y')\in\calS(\Lambda,
    \calC^{\infty,\gamma+\eps}(\bz_y)\pit
    \calC^{\infty,-\gamma+\eps}(\bz_{y'})),$$
\end{definition}
cf. \eqref{not5}; $\pit$ denotes the completed projective tensor product. 
$C^{-\infty}_G(\bz;\Lambda,\gamma)$ is the residual class of the calculus.   
For every choice of $s$, $p$, and $\lambda$, the operator $G(\lambda)$
maps $\hsgpb$ into $\cigb$. For the description of the resolvent we shall
need another class of operator-families. For each fixed $\lambda$, they are
smoothing over $X^\wedge$, yet they have a finite order in $\lambda$: 
\begin{definition}\label{green2}
 Let $\gamma,\mu\in\rz$ and $d>0$. We define 
 $R^{\mu,d}_G(X^\wedge;\Lambda,\gamma)$ as the space of all 
 operator-families $G=G(\lambda)$ that have a kernel with respect to the 
 $\calK^{0,0}_2(X^\wedge)$-scalar product of the form 
  $$k(\lambda,t,x,t',x')=[\lambda]^{\frac{n+1}{d}}
    \tilde{k}(\lambda,[\lambda]^{\frac{1}{d}}t,x,
              [\lambda]^{\frac{1}{d}}t',x'),$$
 where $[\cdot]$ is a smoothed norm-function $($i.e., $[\cdot]$ is 
 smooth, positive on $\cz$ and $[\lambda] = |\lambda|$ for large 
 $\lambda)$ and for some $\eps=\eps(G)>0$
  $$\tilde{k}(\lambda,t,x,t',x')\in S^{\frac{\mu}{d}}_{cl}(\Lambda)\pit
    \calS^{\gamma+\eps}_0(X^\wedge_{(t,x)})\pit
    \calS^{-\gamma+\eps}_0(X^\wedge_{(t',x')}).$$
\end{definition}
In this case, $G(\lambda)$ maps $\calK^{s,\gamma}_p(X^\wedge)$ into 
$\calS^\gamma_0(X^\wedge)$ for any $s$ and $p$. See also the Appendix for 
more information on such operator-families. Trivially, a symbol 
$a\in S^{\frac{\mu}{d}}_{cl}(\Lambda)$ satisfies the estimate 
 \begin{equation}\label{decay}
  |a(\lambda)|\le c\,(1+|\lambda|)^{\frac{\mu}{d}},\qquad
  \lambda\in\Lambda.
 \end{equation}
Recall from Section \ref{diffop} that if $A$ is a cone differential 
operator, then $\lambda-A$ can be written in terms of Mellin symbols taking 
values in the differential operators on $X$, cf.\ \eqref{tildeh}. In that 
case the Mellin symbol is a polynomial in $z$. A general Mellin 
symbol is an entire function with values in 
the pseudodifferential operators on $X$; more precisely:  
\begin{definition}\label{holom}
 For $\mu\in\rz$ and $d>0$ let $M^{\mu,d}_{\calO}(X;\Lambda)$ denote the 
 space of all functions $\tilde{g}(z,\lambda)$, which are holomorphic in 
 $z\in\cz$ with values in $L^{\mu,d}(X;\Lambda)$, and for which
  $$\tilde{g}_\beta(\tau,\lambda):=\tilde{g}(\beta+i\tau,\lambda)\;\in\;
    L^{\mu,d}(X;\rz_\tau\times\Lambda)$$
 is locally bounded as a function of $\beta\in\rz$. This is a Fr\'echet 
 space in a canonical way.
\end{definition}
Let us now state the ellipticity assumptions on $A$, which ensure 
the existence of its resolvent in a keyhole region:
\begin{definition}\label{ellipticity}
 We call $A$ elliptic with respect to the weight $\gamma+\mu$ and the 
 sector $\Lambda_\Delta$, cf.\ \eqref{not4}, if 
 the following two conditions are satisfied: \begin{itemize}
  \item[(E1)] Both the homogeneous principal symbol $\sigma_\psi^\mu(A)$ 
   and the rescaled symbol $\widetilde{\sigma}_\psi^\mu(A)$, cf.\ 
   \eqref{principal2},  have no spectrum in $\Lambda_\Delta$, pointwise on
$T^*({\rm int}\,\bz)\setminus 0$ and $(T^*X\times \rz)\setminus 0$, 
respectively, 
  \item[(E2)] the model cone operator $\widehat{A}$, acting as in \eqref{ahat}, 
   has no spectrum in $\Lambda_\Delta\setminus\{0\}$. 
 \end{itemize}
\end{definition}
If conditions (E1) and (E2) are satisfied, they automatically hold for a 
slightly larger keyhole region (by closedness of the spectrum, compactness 
of $\bz$, and the homogeneity of the rescaled symbol, the homogeneous 
principal symbol, as well as $\widehat{a}(\lambda)$, cf.\ \eqref{twisted}).
Moreover, one can show that condition (E2) implies that \eqref{conormal2} is 
a family of isomorphisms. 

We would like to point out that, although the above conditions seem to 
be quite strong, it follows from more general considerations that they 
are essentially necessary.

Under conditions (E1) and (E2) we can now describe the resolvent of $A$: 
\begin{theorem}\label{parametrix}
 If $A$ is elliptic with respect to $\Lambda_\Delta = \Lambda_{\Delta}(\theta)$
 and $\gamma+\mu$, 
 then $A$ has no spectrum in $\Lambda_\Delta\cap\{|\lambda|>R\}$ for some 
 $R>0$, and for large $\lambda\in\Lambda_\Delta$
  \begin{equation}\label{inverse}
   (\lambda-A)^{-1}=\sigma\left\{t^\mu\opm{\gamma-\frac{n}{2}}(g)(\lambda)+
    G(\lambda)\right\}\sigma_0+(1-\sigma)P(\lambda)(1-\sigma_1)+
    G_\infty(\lambda),
  \end{equation}
 where $\sigma,\sigma_0,\sigma_1\in\cicomp([0,1[)$ are cut-off functions 
 satisfying $\sigma_1\sigma=\sigma_1$, $\sigma\sigma_0=\sigma$, and 
 there exists a $\delta > 0$ such that, for $\Lambda = 
 \Lambda(\delta,\theta)$,
  \begin{itemize}
   \item[i)] 
    $g(t,z,\lambda)=\tilde{g}(t,z,t^\mu\lambda)$ with 
    $\tilde{g}\in \ci(\rpbar,M^{\mu,d}_{\calO}(X;\Lambda))$, 
   \item[ii)] $P(\lambda)\in L^{-\mu,\mu}(\intb;\Lambda)$, cf. 
   \eqref{not6}, 
   \item[iii)] $G(\lambda)\in 
    R^{-\mu,\mu}_G(X^\wedge;\Lambda,\gamma)$, and 
    $G_\infty(\lambda)\in C^{-\infty}_G(\bz;\Lambda,\gamma)$. 
  \end{itemize}
\end{theorem}
In view of the fact that $A$ has compact resolvent (recall that the 
embeddings $\hsgpb\hookrightarrow\calH^{r,\varrho}_p(\bz)$ are compact 
provided $s>r$ and $\gamma>\varrho$), only finitely many points of the 
spectrum of $A$ will lie in $\Lambda$. Thus, after possibly rotating $A$ a 
little and shrinking the keyhole $\Lambda$, we can assume that 
$A$ has no spectrum in $\Lambda$, except perhaps 0. 

Theorem \ref{parametrix} follows from the parametrix 
construction in the parameter-dependent cone algebra
given in \cite{Gil}, Theorems 3.2, 3.4, cf.\ also \cite{EgSc},
Section 9.3.3, Theorem 6. An important observation we can draw from 
this theorem is a norm estimate of the resolvent: 
\begin{proposition}\label{estimate}
 Under the assumptions of Theorem \text{\rm \ref{parametrix}} there exists 
 a constant $c_p\ge 0$ such that for all sufficiently large
 $\lambda\in\Lambda$ 
  $$\|(\lambda-A)^{-1}\|_{\calL(\calH^{0,\gamma}_p(\bz))}\le 
    c_p\,|\lambda|^{-1}.$$
\end{proposition}
\begin{proof}
 We first reduce the case of arbitrary $\gamma$ to the 
 special case $\gamma=\gamma_p=(n+1)(\mbox{$\frac{1}{2}-\frac{1}{p}$})$. To 
 this end let $b\in\ci(\intb)$ be a positive function such that 
 $b(t,x)=t^\nu$, $\nu=\gamma_p-\gamma$, for all $(t,x)\in[0,1]\times X$. 
 Multiplication by $b$ induces isomorphisms 
 $\hsgpb\to\calH^{s,\gamma_p}_p(\bz)$ with inverse induced by $b^{-1}$. 
 Therefore, 
  $$\|(\lambda-A)^{-1}\|_{\calL(\calH^{0,\gamma}_p(\bz))}\sim
    \|b(\lambda-A)^{-1}b^{-1}\|_{\calL(\calH^{0,\gamma_p}_p(\bz))}.$$ 
 But for large $|\lambda|$, 
  $$b(\lambda-A)^{-1}b^{-1}=\sigma\left\{
    t^\mu\opm{\gamma_p-\frac{n}{2}}(T^{-\nu}g)(\lambda)+
    t^\nu G(\lambda)t^{-\nu}\right\}\sigma_0+
    (1-\sigma)bP(\lambda)b^{-1}(1-\sigma_1)+bG_\infty(\lambda)b^{-1},$$
 where $T^{-\nu}g(t,z,\lambda)=g(t,z-\nu,\lambda)$. Since $T^{-\nu}g$ and 
 $bP(\lambda)b^{-1}$ are of the same quality as $g$ and $P(\lambda)$, 
 respectively, and $t^\nu G(\lambda)t^{-\nu}\in 
 R^{-\mu,\mu}_G(X^\wedge;\Lambda,\gamma_p)$ and 
 $bG_\infty(\lambda)b^{-1}\in C^{-\infty}_G(\bz;\Lambda,\gamma_p)$, we can 
 assume from the very beginning that $\gamma=\gamma_p$. 

 The term $G_\infty(\lambda)$ certainly behaves in the right way, since it 
 is rapidly decreasing in $\lambda$. Also the term 
 $(1-\sigma)P(\lambda)(1-\sigma_1)$ is good by the standard 
 Calderon-Vaillancourt theorem. The two remaining terms 
 $t^\mu\opm{\gamma}(g)(\lambda)$ and $G(\lambda)$ we shall consider in the 
 spaces 
  $$\calK^{0,\gamma_p}_p(X^\wedge)=L_p(\rz_+\times X,t^ndtdx).$$ 
 If $\kappa_\varrho$ is the group action from \eqref{kappa}, then 
 $\|\kappa_\varrho\|_{\calL(\calK^{0,\gamma_p}_p(X^\wedge))}=
 \varrho^{\gamma_p}$ 
 for all $\varrho>0$. Hence for an arbitrary operator 
 $T\in\calL(\calK^{0,\gamma_p}_p(X^\wedge))$ we have 
  $$\|T\|_{\calL(\calK^{0,\gamma_p}_p(X^\wedge))}= 
    \|\kappa_\varrho^{-1}T\kappa_\varrho\|_{
    \calL(\calK^{0,\gamma_p}_p(X^\wedge))}.$$
 Now let $G(\lambda)$ have a kernel $k(\lambda)$ as described in Definition 
 \ref{green2} (with $\mu$ replaced by $-\mu$ and $d=\mu$). Then the 
 operator-norm of $G(\lambda)$ is the same as that of 
 $\kappa_{[\lambda]^{1/\mu}}^{-1}G(\lambda)\kappa_{[\lambda]^{1/\mu}}$, 
 which has the kernel $\tilde{k}(\lambda,t,x,t',x')$. But this kernel is 
 $O(|\lambda|^{-1})$ in $\lambda$, cf.\ \eqref{decay}. To treat the last 
 term we can pass to local coordinates, i.e.\ we assume $X=\rz^n$ and 
  $$\tilde{g}_{\frac{1}{2}-\gamma_p}(t,x,\tau,\xi,\lambda)\;\in\;
    S^{-\mu,\mu}(\rpbar\times\rz^n\times\rz^{1+n};\Lambda).$$
 By a tensor product argument we can assume that $\tilde{g}$ is 
 independent of $(t,x)$. Conjugating with 
 $\kappa_{[\lambda]^{1/\mu}}$, we have to show that 
  $$\opm{\gamma_p-\frac{n}{2}}(g')(\lambda),\qquad 
    g'_{\frac{1}{2}-\gamma_p}(t,\tau,\xi,\lambda)=t^\mu
    \tilde{g}_{\frac{1}{2}-\gamma_p}(\tau,\xi,t^{\mu}
    \mbox{$\frac{\lambda}{[\lambda]}$}),$$
 is uniformly bounded in $L_p(\rz_+\times\rz^n,t^ndtdx)$ for large 
 $|\lambda|$. Since then $\frac{\lambda}{[\lambda]}$ is bounded away from 
 zero and infinity, a simple calculation shows that 
  $$|(t\partial_t)^l\partial_\tau^k\partial_\xi^\alpha 
    g'_{\frac{1}{2}-\gamma_p}(t,\tau,\xi,\lambda)|\le 
    c_{kl\alpha}(1+|\tau|+|\xi|)^{-k-|\alpha|}$$
 uniformly in $(t,\tau,\xi,\lambda)$, i.e., $g'(\lambda) \in 
   MS^{0}(\rz_+\times\rz^n\times
   \Gamma_{\frac{n+1}{2}-\gamma_{p}}\times\rz^n)$ uniformly in 
   $\lambda$. Then the result follows, see the end of Section 
   \ref{notation}.
 \end{proof}
There are certain relations between $A$ and $g$ from i) respectively 
$P(\lambda)$ from ii), which we are going to study now. 

Let $U\subset\rz^n$ be a coordinate neighborhood for $X$, whose closure is 
contained in another coordinate neighborhood. Condition (E1) ensures that 
the local symbol 
 \begin{equation}\label{qtilde}
  \tilde{q}^{(\mu)}(t,x,\tau,\xi):=
  t^\mu\sigma_\psi^\mu(A)(t,x,t^{-1}\tau,\xi),
 \end{equation} 
cf.\ \eqref{principal}, exists up to $t=0$ and does not have spectrum in 
$\Lambda_\Delta$ for all $(t,x)\in[0,1]\times\overline{U}$ and all 
$(\tau,\xi)\not=0$. 
\begin{lemma}\label{excise}
 There exists a zero excision function $\chi$ on $\rz$, such that for any 
 $\varphi\in\cicomp(U)$ and any $\sigma\in\cicomp([0,1[)$ 
  $$\varphi(x)\sigma(t)\chi(|\tau,\xi|^2+|\lambda|^{\frac{2}{\mu}})
    \big(\lambda-\tilde{q}^{(\mu)}(t,x,\tau,\xi)\big)^{-1}\;\in\;
    S^{-\mu,\mu}(\rpbar\times\rz^n\times\rz^{1+n};\Lambda).$$    
\end{lemma}
\begin{proof}
 Let $a$ denote the symbol in question. For shortness let us write 
 $y=(t,x)$ and $\eta=(\tau,\xi)$. Since the eigenvalues of 
 $\tilde{q}^{(\mu)}(y,\eta)$ are proportional to $|\eta|^\mu$ (uniformly 
 for $y\in[0,1]\times\overline{U}$) and do not lie in $\Lambda_\Delta$, 
 there exists a constant $c>0$ such that 
 $(\lambda-\tilde{q}^{(\mu)}(y,\eta))^{-1}$ is a smooth function in 
  $$\{(y,\eta,\lambda)\in[0,1]\times\overline{U}\times\rz^{1+n}\times\cz\st
    \lambda\in\Lambda_\Delta\text{ or }|\lambda|\le c|\eta|^\mu\}.$$ 
 Thus, if we choose $\chi$ in such a way that 
 $\chi(|\eta|^2+|\lambda|^{\frac{2}{\mu}})$ vanishes for 
 $|\eta|\le(\delta/c)^{\frac{1}{\mu}}$ and $|\lambda|\le\delta$, then $a$ 
 is smooth on $\rpbar\times\rz^n\times\rz^{1+n}\times\Lambda$ for 
 $\Lambda=\Lambda(\delta,\theta)$. To verify that $a$ is a symbol, it 
 suffices to show that 
  \begin{equation}\label{formula1}
   |a(y,\eta,\lambda)|\le c\,
   (1+|\eta|^2+|\lambda|^{\frac{2}{\mu}})^{\frac{-\mu}{2}}
  \end{equation}
 uniformly in $y\in[0,1]\times\overline{U}$ and 
 $(\eta,\lambda)\in\rz^{1+n}\times\Lambda$. Since $a$ is anisotropic 
 homogeneous of order $(-\mu,\mu)$ for large 
 $(\eta,\lambda)\in\rz^{1+n}\times\Lambda_\Delta$, estimate 
 \eqref{formula1} holds on $\rz^{1+n}\times\Lambda_\Delta$. It also holds 
 for $|\lambda|\le\delta$ and $|\eta|$ sufficiently large, since then 
 $|(\lambda-\tilde{q}^{(\mu)}(y,\eta))^{-1}|$ is  
 $O(|\eta|^{-\mu})$ due to the above described behavior of the eigenvalues. 
 For $|\lambda|$ and $|\eta|$ simultaneously small, estimate 
 \eqref{formula1} holds anyway. 
\end{proof}
For every $\beta \in \rz$, we can associate with $\tilde{g}$ from Theorem 
\ref{parametrix}.i) a local symbol 
 $$\tilde{g}_\beta=\tilde{g}_\beta(t,x,\tau,\xi,\lambda)\;\in\;
   S^{-\mu,\mu}(\rpbar\times\rz^n\times\rz^{1+n};\Lambda).$$
   It is a consequence of the above mentioned parametrix construction 
   in the cone calculus that the principal symbol of $\tilde{g}$ is 
   determined by the inverted principal symbol of $\lambda - A$.
With the notation from Lemma \ref{excise} we indeed have 
 \begin{equation}\label{param2}
  \varphi(x)\sigma(t)\left\{\tilde{g}_\beta-   
  \chi(|\tau,\xi|^2+|\lambda|^{\frac{2}{\mu}})
  \big(\lambda-\tilde{q}^{(\mu)}\big)^{-1}\right\}\;\in\;
  S^{-\mu-1,\mu}(\rpbar\times\rz^n\times\rz^{1+n};\Lambda).
 \end{equation}
Similarly, the local symbols of $P(\lambda)$ from Theorem 
\ref{parametrix}.ii) can be approximated modulo $S^{-\mu-1,\mu}$ in terms 
of the inverted local principal symbol of $\lambda - A$. 

\section{Complex powers of cone differential operators}\label{boimpo}
The aim of this section is to show that a cone differential operator $A$ 
satisfying (E1), (E2), also satisfies condition (A3). More precisely we 
shall show: 
\begin{theorem}\label{bip}
 Let $A$ be elliptic with respect to $\gamma+\mu$ and $\Lambda_\Delta$,  
 having no spectrum in the keyhole $\Lambda=\Lambda(\delta,\theta)$, 
 except perhaps 0. Then one can define $A^z$ as in \eqref{atoz} and 
 there exists a constant $c_p\ge0$ such that for all $z\in\hz$ with $|\re 
 z|$ sufficiently small 
  \begin{equation}\label{bound}
   \|A^z\|_{\calL(\calH^{0,\gamma}_p(\bz))}\le c_p\,e^{\theta|\im z|}. 
  \end{equation}
\end{theorem}
Let us first give a short outline of the proof. In view of Theorem 
\ref{parametrix} we can replace in \eqref{atoz} the resolvent 
$(\lambda-A)^{-1}$ by the right-hand side of \eqref{inverse}. Then one 
obtains four integrals (corresponding to the four terms on the right-hand
side of \eqref{inverse}), each of which has to be estimated as in \eqref{bound}. 
For the one associated with $G_\infty(\lambda)$ this is certainly possible, 
since $G_\infty(\lambda)$ is rapidly decreasing in $\lambda$ and therefore
 $$\Big\|\int_\calC\lambda^zG_\infty(\lambda)\,d\lambda
   \Big\|_{\calL(\calH^{0,\gamma}_p(\bz))}\le
   c'_p\,\delta^{\re z} e^{\theta|\im z|}.$$ 
Obviously $\delta^{\re z}$ is uniformly bounded for small $|\re z|$. 

For the integral connected with $(1-\sigma)P(\lambda)(1-\sigma_1)$ one can 
proceed exactly as in the proof of Theorem 1 of \cite{Seel2}, 
since this term is localized away from the boundary, and there, 
$\calH^{0,\gamma}_p(\bz)$ coincides with the usual $L_p$-spaces (note also 
the remark after formula \eqref{param2}). 

Hence it remains to consider the expressions 
 \begin{equation}\label{twoterms}
  \sigma\int_\calC\lambda^zG(\lambda)\,d\lambda\,\sigma_0,\qquad
  \sigma\int_\calC\lambda^zt^\mu\opm{\gamma-\frac{n}{2}}(g)(\lambda)\,
  d\lambda\,\sigma_0.
 \end{equation}  
  
We shall start with the analysis of the first term. Letting 
 $$\calH^{0,\gamma}_p(X^\wedge)=
   L_p(X^\wedge,t^{(\frac{n+1}{2}-\gamma)p}\mbox{$\frac{dt}{t}$}dx)=
   t^\gamma L_p(X^\wedge,t^{\frac{n+1}{2}p}\mbox{$\frac{dt}{t}$}dx),$$ 
it is obvious from Definition \ref{sobolev} that multiplication with 
any cut-off function $\sigma\in\cicomp([0,1[)$ induces continuous operators 
$\calH^{0,\gamma}_p(\bz)\to\calH^{0,\gamma}_p(X^\wedge)$ and 
$\calH^{0,\gamma}_p(X^\wedge)\to\calH^{0,\gamma}_p(\bz)$. 
 Estimating the 
first term in \ref{twoterms} thus reduces to the following proposition: 
\begin{proposition}\label{greenpart}
 Let $G(\lambda)\in R^{-\mu,\mu}_G(X^\wedge;\Lambda,\gamma)$ and 
 $G_z={\displaystyle\int_\calC}\lambda^zG(\lambda)\,d\lambda$ for 
 $z\in\hz$. Then $G_z\in\calL(\calH^{0,\gamma}_p(X^\wedge))$, and there 
 exists a constant $c_p\ge0$ such that for $|\re z|$ sufficiently small 
  $$\|G_z\|_{\calL(\calH^{0,\gamma}_p(X^\wedge))}\le 
    c_p\,e^{\theta|\im z|}.$$
\end{proposition}
\begin{proof}
 By conjugation with $t^\gamma$ we can assume  
 that $\gamma=0$ (cf.\ the proof of Proposition \ref{estimate}).  
 If we split the integral into three terms according to the decomposition of 
 $\calC$ in \eqref{not3}, the integral over $\calC_2$ can be estimated in 
 the desired way, since $\|G(\lambda)\|$ is bounded on $\calC_2$. By 
 symmetry, $\calC_1$ and $\calC_3$ can be treated in the same way. So we 
 shall assume for the rest of the proof that 
  $$\calC(t)=\calC_1(t)=te^{i\theta},\qquad -\infty<t\le1,$$
 (for notational convenience we replace $\delta$ by 1). Also for 
 convenience we suppress the $x$-variables from the 
 notation. We shall frequently make use of the fact that, substituting   
 $\lambda=\varrho^\mu e^{i\theta}$, we have 
  $$\int_\calC f(\lambda)\,d\lambda=\mu e^{i\theta}\int_1^\infty 
   f(\varrho^\mu e^{i\theta})\varrho^{\mu-1}\,d\varrho.$$ 
 According to Definition \ref{green2}, $G(\lambda)$ for $|\lambda|\ge 1$ is 
 an integral operator (with respect to the scalar product in 
 $\calH^{0,0}_2(X^\wedge)$) with kernel 
  $$k(\lambda,t,s)=|\lambda|^{\frac{n+1}{\mu}}\tilde{k}(\lambda,
    |\lambda|^{\frac{1}{\mu}}t,|\lambda|^{\frac{1}{\mu}}s),$$ 
 where, for some $\varepsilon>0$, 
  $$\tilde{k}(\lambda,t,s)\in S^{-1}(\Lambda) \pit 
    \calS^\eps_0(X^\wedge)\pit \calS^\eps_0(X^\wedge)$$
 (the fact that $\tilde{k}$ is classical will  
 not play a role for the following calculations). Then $G_z$ is an 
 integral operator with kernel 
  \begin{equation}\label{kernelgz}
   k_z(t,s)=\int_\calC\lambda^zk(\lambda,t,s)\,d\lambda.
  \end{equation} 
 Writing $\tilde{k}(\lambda,t,s)=
 (\tilde{\chi}(t)+(1-\tilde{\chi})(t))\tilde{k}(\lambda,t,s)
 (\tilde{\chi}(s)+(1-\tilde{\chi})(s))$ with the 
 characteristic function $\tilde{\chi}$ of $[0,1]$, the proposition will be 
 true, if we can show that in any of the four cases 
 \begin{align}
  k_z(t,s)&=e^{-\theta|\im z|}\int_\calC\lambda^z
            \tilde{\chi}(|\lambda|^{\frac{1}{\mu}}t)k(\lambda,t,s)
            \tilde{\chi}(|\lambda|^{\frac{1}{\mu}}s)\,d\lambda
            \label{kernel1}\\
  k_z(t,s)&=e^{-\theta|\im z|}\int_\calC\lambda^z
            \tilde{\chi}(|\lambda|^{\frac{1}{\mu}}t)k(\lambda,t,s)
            (1-\tilde{\chi})(|\lambda|^{\frac{1}{\mu}}s)\,d\lambda
            \label{kernel2a}\\
  k_z(t,s)&=e^{-\theta|\im z|}\int_\calC\lambda^z
            (1-\tilde{\chi})(|\lambda|^{\frac{1}{\mu}}t)
            k(\lambda,t,s)
            \tilde{\chi}(|\lambda|^{\frac{1}{\mu}}s)\,d\lambda
            \label{kernel2b}\\
  k_z(t,s)&=e^{-\theta|\im z|}\int_\calC\lambda^z
            (1-\tilde{\chi})(|\lambda|^{\frac{1}{\mu}}t)
            k(\lambda,t,s)
            (1-\tilde{\chi})(|\lambda|^{\frac{1}{\mu}}s)\,d\lambda
            \label{kernel3}
 \end{align}
 the associated integral operators are bounded in 
 $\calH^{0,0}_p(X^\wedge)$, uniformly in $-\alpha\le\re z<0$ for some 
 $\alpha>0$. The cases \eqref{kernel2a} and \eqref{kernel2b} are equivalent 
 by symmetry (i.e.\ passing to the adjoint). The proofs of all cases 
 \eqref{kernel1}, \eqref{kernel2a}, and \eqref{kernel3} rely on the 
 following Hardy inequalities: 
  \begin{align}
   \int_0^\infty\Big(\int_0^t g(s)\,ds\Big)^pt^{-1-r}\,dt&\le 
      \mbox{$\left(\frac{p}{r}\right)^p$}
      \int_0^\infty g(t)^pt^{p-1-r}\,dt\label{hardy1}\\
   \int_0^\infty\Big(\int_t^\infty g(s)\,ds\Big)^pt^{-1+r}\,dt&\le 
      \mbox{$\left(\frac{p}{r}\right)^p$}
      \int_0^\infty g(t)^pt^{p-1+r}\,dt\label{hardy2}
  \end{align}
 for any non-negative function $g$ on $\rz_+$ and $r>0$ (cf.\ \cite{StWe}, 
 Lemma 3.14, page 196). To begin with case \eqref{kernel1} we use the fact 
 that, for some fixed $\epsilon > 0$, 
  $$|\tilde{k}(\lambda,t,s)|\le c|\lambda|^{-1}
    t^{-\frac{n+1}{2}+\eps}s^{-\frac{n+1}{2}+\eps}$$ 
 uniformly in $\lambda\in\calC$ and $t,s>0$, to obtain 
  \begin{align*}
   |k_z(t,s)|&\le c\,t^{-\frac{n+1}{2}+\eps}s^{-\frac{n+1}{2}+\eps}
       \int_1^\infty\varrho^{\mu\re z-1+2\eps}
       \tilde{\chi}(\varrho t)\tilde{\chi}(\varrho s)\,d\varrho\\
      &=\mbox{$\frac{c}{\mu\re z+2\eps}$}\left(
       \min(\mbox{$\frac{1}{t},\frac{1}{s}$})^{\mu\re z+2\eps}-1\right)
       \tilde{\chi}(t)\tilde{\chi}(s)
       t^{-\frac{n+1}{2}+\eps}s^{-\frac{n+1}{2}+\eps}.
  \end{align*}
 Since $\mu\re z$ is negative, the factor 
 $\min(\frac{1}{t},\frac{1}{s})^{\mu\re z}$ is uniformly bounded by 1 for 
 $0<s,t\le 1$. If $-\frac{\eps}{\mu}=-\alpha\le\re z<0$ the factor 
 $\frac{c}{\mu\re z+2\eps}$ can be estimated from above by a constant 
 uniformly in $0<s,t\le1$. Since furthermore the kernel function 
 $\tilde{\chi}(t)\tilde{\chi}(s)t^{-\frac{n+1}{2}+\eps}s^{-\frac{n+1}{2}+\eps}$ 
 belongs to $\calH^{0,0}_p(X^\wedge)\otimes\calH^{0,0}_{p'}(X^\wedge)$ 
 and thus induces a continuous operator in $\calH^{0,0}_p(X^\wedge)$, 
 it remains to consider the kernel 
 $t^{-\frac{n+1}{2}+\eps}s^{-\frac{n+1}{2}+\eps}
 \min(\frac{1}{t},\frac{1}{s})^{2\eps}$. Because this kernel is symmetric 
 in $s$ and $t$, indeed it suffices to treat 
  $$k(t,s)=\begin{cases}
            t^{-\frac{n+1}{2}-\eps}s^{-\frac{n+1}{2}+\eps}&:s\le t\\
            0&:s>t
           \end{cases}.$$
 If $G$ denotes the associated integral operator, then 
  \begin{align*}
   \|Gu\|^p_{\calH^{0,0}_p(X^\wedge)}&\le
     \int_0^\infty\Big(\int_0^\infty 
     k(t,s)|u(s)|\,s^nds\Big)^pt^{\frac{n+1}{2}p-1}\,dt
     =\int_0^\infty\Big(\int_0^t s^{\frac{n-1}{2}+\eps}|u(s)|\,ds\Big)^p
     t^{-1-p\eps}\,dt\\
    &\le\mbox{$\left(\frac{p}{p\eps}\right)^p$}\int_0^\infty 
     |u(t)|^p\,t^{\frac{n+1}{2}p-1}dt=
     \mbox{$\left(\frac{1}{\eps}\right)^p$}
     \|u\|^p_{\calH^{0,0}_p(X^\wedge)}
  \end{align*}
 by Hardy's inequality \eqref{hardy1}. This finishes case \eqref{kernel1}. 
 For case \eqref{kernel2a} observe that 
  $$|\tilde{k}(\lambda,t,s)(1-\tilde{\chi})(s)|\le c_N |\lambda|^{-1}
    t^{-\frac{n+1}{2}+\eps}s^{-N}$$ 
 for any $N\in\nz$ uniformly in $\lambda\in\calC$ and $s,t>0$. Then 
  $$|k_z(t,s)|\le c_N t^{-\frac{n+1}{2}+\eps}s^{-N}
    \int_1^\infty\varrho^{\mu\re z+\frac{n-1}{2}+\eps-N}
    \tilde{\chi}(\varrho t)(1-\tilde{\chi})(\varrho s)\,d\varrho.$$
 This expression equals zero if $s\le t$ and for $s>t$ we can estimate 
  \begin{align*}
   |k_z(t,s)|&\le c_N t^{-\frac{n+1}{2}+\eps}s^{-N}\tilde{\chi}(t)
      \int_{1/s}^{1/t}\varrho^{\mu\re z+\frac{n-1}{2}+\eps-N}\,d\varrho
   =\mbox{$\frac{c_N}{\mu\re z+\frac{n+1}{2}+\eps-N}$}
      (k^1_z(t,s)-k^2_z(t,s))
  \end{align*}
 with kernel functions $k^1_z$ and $k^2_z$ given by 
  \begin{align*}
   k^1_z(t,s)&=\tilde{\chi}(t)
               \begin{cases}
                0&:s\le t\\
                \left(\frac{t}{s}\right)^N t^{-n-1-\mu\re z}&:s>t
               \end{cases},\qquad
   k^2_z(t,s)&=\tilde{\chi}(t)
               \begin{cases}
                0&:s\le t\\
                (ts)^{-\frac{n+1}{2}}\left(\frac{t}{s}\right)^\eps
                s^{-\mu\re z}&:s>t
               \end{cases}.
  \end{align*}
  In order to check the uniform boundedness of the integral operator
  $K^{1}_{z}$ associated with $k^{1}_{z}$, $-\alpha \le \re z < 0$, on 
  $\calH^{0,0}_p(X^\wedge)$ we observe that
  $$
  \calH^{0,0}_p(X^\wedge) = t^\beta L_{p}(X^\wedge,t^ndtdx)
  $$
  with $\beta = p(n+1)(\frac{1}{2}-\frac{1}{p})$. The boundedness of 
  $K^{1}_{z}$ is equivalent to the boundedness of $t^\beta K^{1}_{z} 
  t^{-\beta}$ on $L_{p}(X^\wedge,t^ndtdx)$. To show it, we employ 
  Schur's lemma: if $N$ is sufficiently large, then
   $$\int_0^\infty t^\beta k^1_z(t,s)\,s^{-\beta}\,s^nds=
    \tilde{\chi}(t)\,t^{\beta+N-n-1-\mu\re z}\int_t^\infty 
    s^{-N+n-\beta}\,ds=
    \mbox{$\frac{1}{N-n+\beta-1}$} t^{-\mu\re z} \tilde{\chi}(t) \le 1$$
  and 
  \begin{align*} 
   \int_0^\infty t^\beta k^1_z(t,s)s^{-\beta}\,t^ndt=
    s^{-N-\beta}\int_0^{\min(1,s)} t^{\beta+N-1-\mu\re z}\,dt
    =\mbox{$\frac{s^{-N-\beta}}{N-\mu\re z+\beta}$}
    \min(1,s)^{\beta+N-\mu\re z}\le 1.
  \end{align*} 
  To handle $k_z^2$ first observe that we can drop the factor $s^{-\mu\re 
 z}$, since this is uniformly bounded by 1 in $s\le 1$ and $\re z<0$, and 
 if $s\ge 1$ and $-\frac{\eps}{2\mu}\le\re z<0$ then 
 $\left(\frac{t}{s}\right)^{\eps/2}s^{-\mu\re z}\le 1$ for $0\le t\le 1$ 
 (for $t\ge1$ anyway $k_z^2(t,s)=0$). Thus we can assume that
  $$k^2_z(t,s)=k^2(t,s)=
               \tilde{\chi}(t)
               \begin{cases}
                0&:s\le t\\
                (ts)^{-\frac{n+1}{2}}\left(\frac{t}{s}\right)^\eps=
                t^{-\frac{n+1}{2}+\eps}s^{-\frac{n+1}{2}-\eps}&:s>t
               \end{cases}.$$
 But then Hardy's inequality \eqref{hardy2} shows that the integral 
 operator associated with the kernel $k^2_{z}$ is continuous in 
 $\calH^{0,0}_p(X^\wedge)$ with operator norm bounded by $\frac{1}{\eps}$. 
 This finishes case \eqref{kernel2a}. For the final case \eqref{kernel3} we 
 use that 
  $$|(1-\tilde{\chi})(t)\tilde{k}(\lambda,t,s)(1-\tilde{\chi})(s)|\le    
    c_N |\lambda|^{-1}t^{-N}s^{-N}$$ 
 for any $N\in\nz$ uniformly in $\lambda\in\calC$ and $s,t>0$. Then 
  \begin{align*}
   |k_z(t,s)|&\le c_N t^{-N}s^{-N}
     \int_1^\infty\varrho^{\mu\re z+n-2N}
     (1-\tilde{\chi})(\varrho t)(1-\tilde{\chi})(\varrho s)\,d\varrho\\
   &=-\mbox{$\frac{c_N}{\mu\re z+n+1-2N}$}t^{-N}s^{-N}
    \max\mbox{$\left(1,\frac{1}{t},\frac{1}{s}\right)$}^{\mu\re z+n+1-2N}.
 \end{align*}
The factor in front is obviously uniformly bounded in $\re z<0$ 
for $N$ sufficiently large. Since $\mu\re z$ is negative, 
 \begin{align*}
  |(1-\tilde{\chi})(t)k_z(t,s)(1-\tilde{\chi})(s)|&\le
    c\,(1-\tilde{\chi})(t)t^{-N}s^{-N}(1-\tilde{\chi})(s),\\
  |\tilde{\chi}(t)k_z(t,s)(1-\tilde{\chi})(s)|&\le
    c\,\tilde{\chi}(t)t^{N-n-1}s^{-N}(1-\tilde{\chi})(s),\\
  |(1-\tilde{\chi})(t)k_z(t,s)\tilde{\chi}(s)|&\le
    c\,(1-\tilde{\chi})(t)t^{-N}s^{N-n-1}\tilde{\chi}(s).
 \end{align*}
 All these kernel functions belong to 
 $\calH^{0,0}_p(X^\wedge)\otimes\calH^{0,0}_{p'}(X^\wedge)$ for 
 sufficiently large $N$ and thus induce continuous operators in 
 $\calH^{0,0}_p(X^\wedge)$. Hence it remains to investigate 
 $\tilde{\chi}(t)k_z(t,s)\tilde{\chi}(s)$ and by symmetry even 
  $$k(t,s)=\begin{cases}
            0&:s\le t\\
            s^{-N}t^{N-n-1}&:s>t
           \end{cases}.$$
 Again Hardy's inequality \eqref{hardy2} shows that the associated 
 operator is $\calH^{0,0}(X^\wedge)$-continuous. 
\end{proof}
 
    We consider now the second term in (\ref{twoterms}).
    Using a partition of unity on $X$ with any two functions supported 
    in a single coordinate neighborhood, we can assume $X = \rz^n$ and use local 
    symbols compactly supported in $x$. To complete the proof of theorem
    (\ref{bip}) we make use of the
    decomposition (\ref{param2}) and of the fact (see Section \ref{notation}) 
    that operators defined by means of symbols 
    $a \in MS^{0}(\rz_+\times\rz^n\times
    \Gamma_{\frac{n+1}{2}-\gamma}\times\rz^n)$ are bounded on
    $\calH^{0,\gamma}_p(\rz_+\times\rz^n)$, with norm estimated in 
    terms of the seminorms associated to $a$. We treat the
    homogeneous principal symbol of $g$ and the lower order part separately. 
    
\begin{lemma}\label{lower}
 Let $\tilde{b}\in 
 S^{-\mu-1,\mu}(\overline{\rz}_+\times\rz^n\times
 \Gamma_{\frac{n+1}{2}-\gamma}\times\rz^{n};\Lambda)$ 
 be compactly supported in $t$ and
  $$b_z(t,x,\mbox{$\frac{n+1}{2}$}-\gamma+i\tau,\xi)=
  t^\mu\int_\calC\lambda^z
    \tilde{b}(t,x,\mbox{$\frac{n+1}{2}$}-\gamma+i\tau,\xi,t^\mu\lambda)\,
    d\lambda.$$
    For $\re z<0$ this defines a symbol $b_z\in 
    MS^{0}(\rz_+\times\rz^n\times
    \Gamma_{\frac{n+1}{2}-\gamma}\times\rz^n)$, and the symbol estimates of 
 $e^{-\theta|\im z|}b_z$ are uniform in $-1\le\re z<0$. Consequently, 
  $$\|\text{\rm op}_M^{\gamma-\frac{n}{2}}(b_z)\|_{\mathcal{L}
     (\calH^{0,\gamma}_p(\rz_+\times\rz^n))}\le 
    c_{p} \, e^{\theta|\im z|}$$ 
 uniformly in $-1\le\re z<0$. 
\end{lemma}
\begin{proof}
     Without loss of generality, we can set $\gamma = \frac{n+1}{2}$. We 
     have to show that 
     \begin{equation}
       \label{ms0estimate}
           | \partial^l_{\tau} (t\partial_{t})^k 
       \partial^\alpha_{\xi} \partial^\beta_{x}
       b_z(t,x,i\tau,\xi) | 
       e^{-\theta|\im z|}\langle\tau,\xi \rangle^{l+|\alpha|}
     \end{equation}
     is uniformly bounded for $t>0$, $x \in \rz^n$, $\tau\in\rz$ and
     $-1\le\re z<0$. The totally characteristic derivatives in $t$ 
     can be handled very simply, observing that $t\partial_t t^\mu=\mu t^\mu$,
  $$t\partial_t\left(\tilde{b}
    (t,x,i\tau,\xi,t^\mu\lambda)\right)=
    (t\partial_t\tilde{b})(t,x,i\tau,\xi,t^\mu\lambda)+
    \mu(\lambda\partial_\lambda\tilde{b})
    (t,x,i\tau,\xi,t^\mu\lambda)$$
  and both symbols $t\partial_t\tilde{b}$ and 
 $\lambda\partial_\lambda\tilde{b}$ are of the same type as $\tilde{b}$.
 Since the derivatives with respect to $x$, $\tau$ and $\xi$ can be 
 taken under the integral sign, it suffices to assume $\tilde{b}\in 
 S^{-\mu-1-k,\mu}(\overline{\rz}_+\times\rz^n\times\Gamma_{0}\times\rz^{n};
 \Lambda)$ and to show that 
  $$|b_z(t,x,i\tau,\xi)|\le c \, e^{\theta|\im 
    z|}\langle\tau, \xi \rangle^{-k} $$
 uniformly in $t>0$, $x \in \rz^n$, $\tau\in\rz$ and $-1\le\re z<0$.
 By hypothesis, we have
  $$|b_z(t,x,i\tau,\xi)|\le 
    c \, t^\mu\int_\calC|\lambda^z|
    (1+\tau^2+|\xi|^{2}+|t^\mu\lambda|^{2/\mu})^{(-\mu-1-k)/2}\,d\lambda,$$
    and on $\calC$ we can estimate $|\lambda^z|$ from above by 
   $\delta^{\re z} e^{\theta|\im z|}$.   
 The transformation $\varrho=t^\mu\lambda$ yields
  $$|b_z(t,x,i\tau,\xi)|\le 
    c\,\delta^{\re z} e^{\theta|\im z|}\langle\tau,\xi\rangle^{-k}
    \int_{t^\mu\calC}(1+|\varrho|)^{-1-\frac{1}{\mu}}\,d\varrho,$$
 where $t^\mu\calC$ means the path $\calC(t^\mu\delta,\theta)$. Since 
 the  support of $\tilde{b}$ is compact, we may assume without loss of generality 
 that $0<t\le 1$. Then we obtain the estimate 
  $$\int_{t^\mu\calC}(1+|\varrho|)^{-1-\frac{1}{\mu}}\,d\varrho\le 
    2\pi\delta+\int_{\calC_\Delta}(1+|\varrho|)^{-1-\frac{1}{\mu}}\,d\varrho<
    +\infty,$$
 and the statement follows, since 
 $\delta^{\re z}$ is uniformly bounded in $-1\le\re z<0$. 
\end{proof}

  \begin{proposition}
 Let $\tilde{g}=\tilde{g}(t,x,z,\xi,\lambda)$ be a local symbol 
 associated to the Mellin symbol $\tilde{g}$ of $(\lambda-A)^{-1}$ of Theorem 
 $\ref{parametrix}.i)$ and let
  $$g_z(t,x,\mbox{$\frac{n+1}{2}-\gamma$}+i\tau,\xi)=
  \sigma(t)t^\mu\int_\calC\lambda^z
    \tilde{g}(t,x,\mbox{$\frac{n+1}{2}-\gamma$}+i\tau,\xi,t^\mu\lambda)\,d\lambda$$
 with some cut-off function 
 $\sigma\in\mathcal{C}^\infty_{\textrm{comp}}([0,1[)$. For 
 $\re z<0$ this defines a symbol 
 $g_z\in MS^{0}(\rz_+\times\rz^n\times
    \Gamma_{\frac{n+1}{2}-\gamma}\times\rz^n)$, and the
 symbol estimates 
 of $e^{-\theta|\im z|}g_z$ are uniform in $-1\le\re z<0$. In particular, 
  $$\|\text{\rm op}_M^{\gamma-\frac{n}{2}}(g_z)\|_{\mathcal{L}
  (\calH^{0,\gamma}_p(\rz_+\times\rz^n))}\le 
    c_{p} \, e^{\theta|\im z|}$$ 
 uniformly in $-1<\re z<0$. 
\end{proposition}  
\begin{proof}
    Without loss of generality let $\gamma = 
    \frac{n+1}{2}$. We shall also suppress $\sigma$ from the notation and instead 
 assume that $0<t\le1$. We can also assume $x$ confined to a compact
 subset of $\rz^n$. By (\ref{param2}),
  $$\tilde{g}(t,x,i\tau,\xi,\lambda)=
    \chi(|\tau,\xi|^{2} + |\lambda|^\frac{2}{\mu})
    (\lambda-\tilde{q}^{(\mu)}(t,x,i\tau,\xi))^{-1}
    \quad\text{mod}\quad 
    S^{-\mu-1,\mu}(\overline{\rz}_+\times\rz^n\times\rz^{1+n};\Lambda),$$
 where $\tilde{q}^{(\mu)}$ denotes a local symbol of $A$ as defined in 
 (\ref{qtilde}). In view of Lemma \ref{lower}, we therefore may assume that 
  \begin{align*}
   g_z(t,x,i\tau,\xi)&=t^\mu\int_\calC\lambda^z
    \chi(|\tau,\xi|^{2}+t^2|\lambda|^\frac{2}{\mu})
    (t^\mu\lambda-\tilde{q}^{(\mu)}(t,x,i\tau,\xi))^{-1}\,d\lambda\\
   &=t^{-\mu 
   z}\int_{t^\mu\calC}\varrho^z\chi(|\tau,\xi|^{2}+|\varrho|^\frac{2}{\mu})
    (\varrho-\tilde{q}^{(\mu)}(t,x,i\tau,\xi))^{-1}\,d\varrho,
  \end{align*}
 where we have used the substitution $\varrho=t^\mu\lambda$. We have 
 to estimate this expression as in (\ref{ms0estimate}).
 The factor $t^{-\mu z}$ behaves correctly, since
 $(t\partial_t)^k t^{-\mu z}=(-\mu z)^k t^{-\mu z}$ is uniformly bounded
 in $0<t\le1$ and $-1\le\re z<0$. For $(\tau,\xi) \not=0$ we have 
  $$\text{spec}(\tilde{q}^{(\mu)}(t,x,i\tau,\xi))\subset 
    \{\lambda\in\cz\st c_1|\tau,\xi|^\mu\le|\lambda|\le 
    c_2|\tau,\xi|^\mu
    \text{ and }|\text{\rm arg}\,\lambda|<\theta\}$$
 with suitable constants $c_1$ and $c_2$. Thus for large enough 
 $|\tau,\xi|$ we have $\chi(|\tau,\xi|^{2}+t^{-2\mu}|\varrho|^{2})=1$ 
 and the spectrum of $\tilde{q}^{(\mu)}(t,x,i\tau,\xi)$ is located to the 
 right of the path $\calC$. 
 By Cauchy's theorem we can then replace the path $t^\mu\calC$ by $\calC$, 
 and obtain for large $|\tau,\xi|$
  $$\int_\calC\varrho^z
    (\varrho-\tilde{q}^{(\mu)}(t,x,i\tau,\xi))^{-1}\,d\varrho=
    2\pi i \, \tilde{q}^{(\mu)}(t,x,i\tau,\xi)^z.$$ 
 Then we can estimate (as in \cite{Seel2}, (2.9))
  $$|\partial_\tau^l(t\partial_t)^k 
     \partial^\alpha_{\xi} \partial^\beta_{x}
     g_z(t,x,i\tau,\xi)|\le 
    p(|z|)e^{\theta|\im z|}\langle\tau,\xi\rangle^{\mu\re z-l-|\alpha|}\le 
    p(|z|)e^{\theta|\im z|}\langle\tau,\xi\rangle^{-l-|\alpha|}$$ 
 with a polynomial $p$. However, since we can replace $\theta$ by 
 $\theta-\varepsilon$ for some 
 $\varepsilon>0$ (as noted in the 
 comments on conditions (E1) and (E2)), this yields the uniform 
 symbol estimates of $g_z$ for large $|\tau,\xi|$. 
 
 For small $|\tau,\xi|$, we now shall show that 
 \begin{equation}\label{gz1}
     g_{z}(t,x,i\tau,\xi)=
     t^{-\mu z}\int_{\Upsilon(t)}\varrho^za(t,x,\tau,\xi;\varrho)\,d\varrho,
 \end{equation} 
 where we have set
 $$
       a(t,x,\tau,\xi;\varrho) = \chi(|\tau,\xi|^{2}+|\varrho|^\frac{2}{\mu})
    (\varrho-\tilde{q}^{(\mu)}(t,x,i\tau,\xi))^{-1} \in 
    S^{-\mu,\mu}(\rpbar\times\rz^{n}\times\rz^{1+n};\Lambda)
    $$
 and $\Upsilon(t)$ is the path given in the following picture (with 
 $r_{0} > 0$ to be chosen appropriately):
   
 
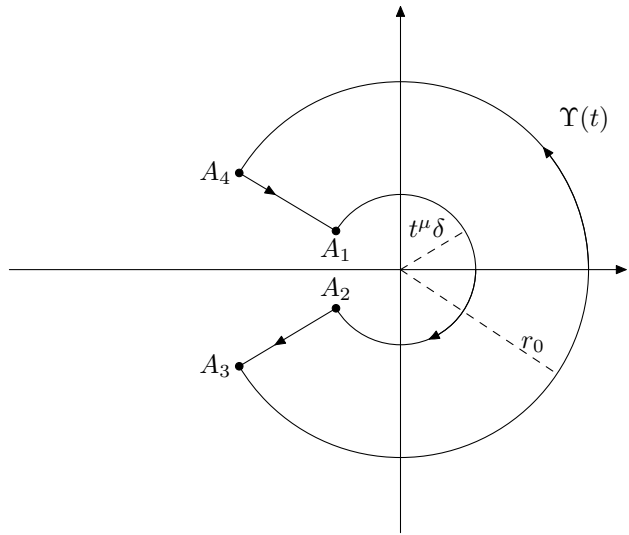
\begin{center}
     \includegraphics{figures.3}
 \end{center}
 
 In fact, the difference of both sides from \eqref{gz1} equals
 \begin{equation}\label{gz2}
     \alpha(r) = \int_{\calC(r,\theta)}\varrho^z
       a(t,x,\tau,\xi;\varrho)
       \,d\varrho \quad \mbox{ for $r = r_{0}$.}
 \end{equation}
 Since, for small $|\tau, \xi|$, the spectrum of 
 $\tilde{q}^{(\mu)}(t,x,i\tau,\xi)$ is contained in some ball of 
 finite radius, $a(t,x,\tau,\xi;\varrho)$ is holomorphic in $\varrho$ 
 for $|\varrho| \ge r_{0}$, if $r_{0}$ is chosen large enough. Thus 
 $\alpha(r) = \alpha(r_{0})$ for all $r \ge r_{0}$, by Cauchy's theorem.
 For any fixed $z$ and $(t,x,\tau,\xi)$ the integrand in \eqref{gz2}
 is $O(|\varrho|^{-1+ \re z})$ for $|\varrho|\to\infty$ and,
 on the radial part of $\calC(r)$, the integrand is 
 $O(r^{\re z})$. Hence, $\alpha(r_{0}) = \lim\limits_{r\to +\infty} 
 \alpha(r) = 0$, and \eqref{gz1} holds.
 
 To estimate the right-hand side of \eqref{gz1},
 we split the integral into four parts, 
 which we briefly analyze separately. 
 First of all, observe that  
 $|\varrho^z|$ can be estimated from above by $e^{\theta|\im 
 z|}(t^\mu\delta)^{\re z}$ on the whole path. This and the fact that 
 $a(t,x,\tau,\xi;\varrho) \in 
 S^{-\mu,\mu}(\rpbar\times\rz^{n}\times\rz^{1+n};\Lambda)$ are enough to get
 the desired  estimates for the terms obtained integrating along the two arcs
 $\stackrel{\frown}{A_{1}A_{2}}$ and $\stackrel{\frown}{A_{3}A_{4}}$, 
 since they can be treated with essentially the same technique we used 
 to prove Lemma \ref{lower}. The term obtained integrating along 
 $\overline{A_{2}A_{3}}$ is 
 $$
       b(t,x,\tau,\xi) =
        \int_{t^\mu \delta}^{r_{0}} (s e^{i\theta})^z
       a(t,x,\tau,\xi;s e^{i\theta})
       e^{i\theta}\,ds.
 $$

 The derivatives with respect to $x$, $\xi$ and $\tau$ can be taken 
 under the integral sign, so that we could again start with a symbol 
 $a \in S^{-\mu-k,\mu}(\rpbar\times\rz^{n}\times\rz^{1+n};\Lambda)$ and prove 
 that, for any $l$,
 $$
 | (t \partial_{t})^l b(t,x,\tau,\xi) | \le c_{l} e^{\theta |\im z|} 
 \langle \tau,\xi \rangle^{-k}
 $$
 uniformly in $-1 < \re z < 0$. This is true for $l = 0$, as one can 
 easily check. For $l=1$ we get
 $$
 t\partial_{t} b(t,x,\tau,\xi) = 
 \int_{t^\mu \delta}^{r_{0}} (s e^{i\theta})^z
       (t\partial_{t}a)(t,x,\tau,\xi;s e^{i\theta})
       e^{i\theta}\,ds -
       \mu (t^\mu \delta e^{i\theta})^{z+1}
       a(t,x,\tau,\xi;t^\mu \delta e^{i\theta}),
 $$
 and this also satisfies the desired estimate. In fact, the first term is 
 of the same kind as $b$, while for the second it suffices to use the 
 definition of $S^{-\mu-k,\mu}(\rpbar\times\rz^{n}\times\rz^{1+n};\Lambda)$
 and the fact that $0<t\le1$.
 The result for arbitrary $l$ can be proved by induction and, obviously, 
 the contribution obtained integrating along $\overline{A_{4}A_{1}}$ behaves 
 in a completely similar way. This yields the desired symbol estimates for
 small $|\tau,\xi|$ and finishes the proof. 
\end{proof} 
\begin{remark}
 Let us point out that the proof of Theorem \text{\rm\ref{bip}} only makes 
 use of 
 assumption \text{\rm(E1)}, the structure of the resolvent 
 \eqref{parametrix}, and the fact that $\text{\rm 
 spec}\,A\cap\Lambda=\emptyset$. It does not use that the domain $\calD(A)$ 
 equals $\calH^{\mu,\gamma+\mu}_p(\bz)$. Therefore, Theorem 
 \text{\rm\ref{bip}} holds 
 true for $A$ considered on other domains, as long as \text{\rm(E1)}, 
 \eqref{parametrix}, and $\text{\rm spec}\,A\cap\Lambda=\emptyset$ are 
 satisfied. 
\end{remark}

Complex powers of Fuchs-type differential operators have been studied 
recently also by Loya \cite{Loya}. He applies Melrose's $b$-calculus 
and focuses on the analytic properties of the kernels in the spirit of 
Seeley \cite{Seel0}.

In the next section we shall investigate the possible closed extensions of 
$A$, and use the previous remark to obtain an analogue of Theorem 
\ref{bip} for the maximal extension of $A$. 

\section{Closed extensions of cone differential operators}\label{closed} 
Let $A$ be a cone differential operator, which is elliptic with respect to 
$\gamma+\mu$ in the sense of Remark \ref{coneell}. If we consider $A$ as 
the unbounded operator in $\calH^{0,\gamma}_p(\bz)$ with domain 
$\cicomp(\intb)$, its closure $\amin=\amin^{\gamma,p}$ is given by 
 $$\calD(\amin)=\calH^{\mu,\gamma+\mu}_p(\bz),$$
and the maximal closed extension $\amax=\amax^{\gamma,p}$ by 
 $$\calD(\amax)=\left\{u\in\calH^{0,\gamma}_p(\bz)\st
   Au\in\calH^{0,\gamma}_p(\bz)\right\}.$$
Note that in \eqref{extension} we simply wrote $A$ instead of $\amin$. 
Taking into account the duality of $\calH^{0,\gamma}_p(\bz)$ and 
$\calH^{0,-\gamma}_{p'}(\bz)$, the following lemma is valid: 
\begin{lemma}\label{adjoint}
 If $A^t$ is the formal adjoint of $A$ with respect to the scalar product 
 of $\calH^{0,0}_2(\bz)$, then 
  $$(\amin^{\gamma,p})^*=(A^t)_{\max}^{-\gamma,p'},\qquad
    (\amax^{\gamma,p})^*=(A^t)_{\min}^{-\gamma,p'}.$$
 We shall write this more shortly as $\amin^*=\amax^t$ and 
 $\amax^*=\amin^t$. 
\end{lemma}
A proof of the above statements in case $p=2$ and $\gamma=0$ is given in 
\cite{Lesc}. The argument in the general case is analogous. As a simple 
consequence, 
 $$(\lambda-\amax)^{-1}=[(\overline{\lambda}-\amin^t)^{-1}]^*$$
whenever one of both sides exists. Since the structure of the resolvent of 
$A=\amin$ as given in Theorem \ref{parametrix} is invariant under passing 
to the adjoint, we obtain the following theorem:
\begin{theorem}\label{parametrix2}
 If $\amin^t$ is elliptic with respect to $\Lambda_\Delta$ and 
 $-\gamma+\mu$, then $\amax$ has no spectrum in 
 $\Lambda_\Delta\cap\{|\lambda|>R\}$ for some $R>0$, and for large 
 $\lambda\in\Lambda_\Delta$
  $$
   (\lambda-\amax)^{-1}=\sigma
    \left\{t^\mu\opm{\gamma-\frac{n}{2}}(g)(\lambda)+
           G(\lambda)\right\}\sigma_0+
    (1-\sigma)P(\lambda)(1-\sigma_1)+G_\infty(\lambda),
  $$
  where $\sigma,\sigma_0,\sigma_1\in\cicomp([0,1[)$ are cut-off functions 
 satisfying $\sigma_1\sigma=\sigma_1$, $\sigma\sigma_0=\sigma$, and 
  \begin{itemize}
   \item[i)] 
    $g(t,z,\lambda)=\tilde{g}(t,z,t^\mu\lambda)$ with 
    $\tilde{g}\in \ci(\rpbar,M^{-\mu,\mu}_{\calO}(X;\Lambda))$, 
   \item[ii)] $P(\lambda)\in L^{-\mu,\mu}(\intb;\Lambda)$, 
   \item[iii)] $G(\lambda)\in 
    R^{-\mu,\mu}_G(X^\wedge;\Lambda,\gamma)$, and 
    $G_\infty(\lambda)\in C^{-\infty}_G(\bz;\Lambda,\gamma)$. 
  \end{itemize} 
\end{theorem}
Proceeding exactly as in Proposition \ref{estimate} and Section 
\ref{boimpo}, we can prove a norm estimate for the complex powers of 
$\amax$: 
\begin{theorem}\label{bip2}
 Let $\amin^t$ be elliptic with respect to $-\gamma+\mu$ and 
 $\Lambda_\Delta$, having no spectrum in the keyhole 
 $\Lambda=\Lambda(\delta,\theta)$ except perhaps 0. Then one can define 
 $\amax^z$ as in \eqref{atoz} and there exists a constant $c_p\ge0$ such 
 that for all $z\in\hz$ with $|\re z|$ sufficiently small 
  $$
   \|\amax^z\|_{\calL(\calH^{0,\gamma}_p(\bz))}\le c_p\,e^{\theta|\im z|}. 
  $$
\end{theorem}
Of course, it is desirable to express the ellipticity assumptions made on 
$\amin^t$ in the previous two theorems purely in terms of $\amax$. This 
can be done as follows.  
\begin{remark}\label{ellipticity2}
 If $\widehat{A}$ is the model cone operator on $X^\wedge$ associated 
 with $A$ and $\widehat{A}_{\max}=\widehat{A}_{\max}^{\gamma,p}$ is the 
 closed operator given by 
  $$\calD(\widehat{A}_{\max})=\left\{u\in\calK^{0,\gamma}_p(X^\wedge)\st
    \widehat{A}u\in\calK^{0,\gamma}_p(X^\wedge)\right\},$$
 then $\amin^t$ is elliptic with respect to $-\gamma+\mu$ and 
 $\Lambda_\Delta$ if and only if $A$ satisfies condition \text{\rm (E1)} 
 and 
  \begin{itemize}
   \item[(E$2'$)] $\widehat{A}_{\max}$ has no spectrum in 
    $\Lambda_\Delta\setminus\{0\}$.
   \end{itemize}
\end{remark}
The previous remark holds true, since similar to Lemma \ref{adjoint}, 
$\widehat{A}_{\min}^*=\widehat{A}_{\max}^t$ and 
$\widehat{A}_{\max}^*=\widehat{A}_{\min}^t$. 

It can be shown that $\calD(\amax)$ differs from $\calD(\amin)$ by a finite 
dimensional space (for the case $p=2$ see \cite{Lesc}), 
 $$\calD(\amax)=\calD(\amin)\oplus V,\qquad\dim V<\infty.$$
More precisely, the dimension of $V$ only depends on the conormal symbol of 
$A$, 
 \begin{equation}\label{gohberg}
   \dim V=\mathop{\mbox{\Large$\sum$}}_{-2<\re z-\frac{n+1}{2}+\gamma<0}
   M(\sigma^\mu_M(A),z),
 \end{equation}
where $M(h,z)$ denotes the multiplicity in $z$ in the sense of \cite{GoSi} 
of a function $h$, which is holomorphic in a punctured neighborhood of 
$z$. Moreover, $V$ consists of smooth functions of the form 
 $$\omega(t)\mathop{\mbox{\Large$\sum$}}_{j=0}^N
   \mathop{\mbox{\Large$\sum$}}_{k=0}^{k_j} 
   c_{jk}(x)t^{-p_j}(\log t)^k,
   \qquad c_{jk}\in\ci(X);$$
the coefficients $c_{jk}$, the exponents $p_j\in\cz$ 
($\frac{n+1}{2}-\gamma-\mu\le\re p_j<\frac{n+1}{2}-\gamma$), and 
$k_j,N\in\nz_0$ are determined by $A$. In particular,
the only closed extensions of $A$ are the operators $A_{W}$ given by 
 $$\calD(A_W)=\calD(\amin)\oplus W,\qquad W\le V.$$ 
In this notation, $\amin=A_{\{0\}}$ and $\amax=A_V$. Correspondingly, 
 $$\calD(\widehat{A}_{\max})=\calD(\widehat{A}_{\min})\oplus\widehat{V},\qquad
   \dim\widehat{V}=\dim V,$$
and all closed extensions $\widehat{A}_{\widehat{W}}$ are given by 
 $$\calD(\widehat{A}_{\widehat{W}})=\calD(\widehat{A}_{\min})\oplus\widehat{W},\qquad
   \widehat{W}\le\widehat{V}.$$
\begin{remark}\label{unique}
 If the conormal symbol $\sigma^\mu_M(A)(z)$ is invertible for all 
 $\frac{n+1}{2}-\gamma-\mu<\re z<\frac{n+1}{2}-\gamma$, then $\dim 
 V=\dim\widehat{V}=0$ by \eqref{gohberg}, and both $A$ and $\widehat{A}$ 
 have only one closed extension in $\calH^{0,\gamma}_p(\bz)$ and 
 $\calK^{0,\gamma}_p(X^\wedge)$, respectively.  
\end{remark}

\section{Example: The Cauchy problem for  Laplacians}\label{laplace} 
Let $g(t)$ be a family of metrics on $X$, depending smoothly on a parameter 
$t\in\rpbar$, and $\Delta_X(t)$ the corresponding Laplacian on $X$. If we 
equip $\intb$ with a metric that coincides with $dt^2+t^2g(t)$ near $t=0$, 
the associated Laplacian $\Delta$ is near the boundary given by 
 $$t^{-2}\left\{(t\partial_t)^2+(n-1+tG^{-1}(t)(\partial_tG)(t))t\partial_t
   +\Delta_X(t)\right\},$$
where $G=(\det(g_{ij}))^{\frac{1}{2}}$ and $n=\dim X$. Hence $\Delta$ is 
a cone differential operator in the sense of \eqref{diffop1}. We shall 
prove the following theorem: 
\begin{theorem}\label{example}
 Let $\Delta$ be the Laplacian on $\intb$ in the above sense, $1<p<\infty$
 such that
  \begin{equation}\label{dimension}
   2\max(p,p')-1<n=\dim\partial\bz. 
  \end{equation} 
 If $\gamma_p=(n+1)(\frac{1}{2}-\frac{1}{p})$, then $\Delta$ defined on 
 $\cicomp(\intb)$ has for any $1<q<\infty$ a unique closed extension 
 $\Delta_{p,q}$ in $\calH^{0,\gamma_p}_q(\bz)$, which is given by
  $$\calD(\Delta_{p,q})=\calH^{2,\gamma_p+2}_q(\bz).$$ 
 Moreover, $-\Delta_{p,q}$ is elliptic with respect to $\gamma_p+2$ and  
 any sector $\Lambda_\Delta\subset\cz\setminus\rz_+$. 
\end{theorem} 
\begin{proof}
 Let us set $A=-\Delta$. The rescaled symbol of $A$ is  
  $$\tilde{\sigma}^{2}_{\psi}(A)(x,\tau,\xi)=\tau^2+|\xi|_x,$$
 where $|\xi|$ refers to the metric $g(0)$ on $X$. Hence $A$ satisfies the 
 ellipticity condition (E1) for any $\Lambda_\Delta$ in question. The 
 conormal symbol of $A$, cf.\ \eqref{conormal} and \eqref{conormal2}, is 
  $$\sigma^2_M(A)(z)=-z^2+(n-1)z-\Delta_X(0):H^s(X)\longrightarrow 
    H^{s-2}(X).$$
 If $0=\lambda_0 \ge \lambda_1 \ge \ldots$ are the eigenvalues of 
 $\Delta_X(0)$, then $\sigma^2_M(A)(z)$ is not bijective if and only if 
  $$z\in\left\{\mbox{$\frac{n-1}{2}\pm
    (\frac{(n-1)^2}{4}-\lambda_j)^{\frac{1}{2}}$}\st j\in\nz_0\right\}.$$
 Note that, in particular, $\sigma^2_M(A)(z)$ is invertible for all $z$ 
 with $0<\re z<n-1$, and thus by condition \eqref{dimension} for all $z$ 
 with $\frac{n+1}{2}-\gamma_p-2\le\re z\le\frac{n+1}{2}-\gamma_p$. This 
 shows that $A$ is elliptic with respect to $\gamma_p+2$ in the sense of 
 Remark \ref{coneell} and has only one closed extension 
  $$A_{p,q}:\calH^{2,\gamma_p+2}_q(\bz)\subset\calH^{0,\gamma_p}_q(\bz)
    \longrightarrow\calH^{0,\gamma_p}_q(\bz)$$ 
 by Remark \ref{unique}. The model cone operator is 
  $$\widehat{A}=-t^{-2}\left\{(t\partial_t)^2+(n-1)(t\partial_t)+
    \Delta_X(0)\right\},$$
 i.e.\ $-\widehat{A}$ is the Laplacian on $X^\wedge$ with respect to the 
 metric $dt^2+t^2g(0)$. As before, $\widehat{A}$ has a unique closed 
 extension 
  $$\widehat{A}_{p,q}:\calK^{2,\gamma_p+2}_q(X^\wedge)\subset
    \calK^{0,\gamma_p}_q(X^\wedge)  
    \longrightarrow\calK^{0,\gamma_p}_q(X^\wedge).$$ 
 Since $\widehat{A}$ is symmetric and non-negative, $\widehat{A}_{2,2}$ is 
 self-adjoint and $\text{\rm spec}(\widehat{A}_{2,2})\subset\rpbar$. Let us 
 show that 
  $$\text{\rm spec}(\widehat{A}_{p,q})\subset\rpbar\qquad
    \forall\;1<q<\infty, \mbox{ $p$ satisfying \eqref{dimension}}.$$
 By Corollary 3.15 of \cite{ScSe1} (in the version for 
 operators in the cone algebra $C^\mu(X^\wedge;(\gamma,\gamma-\mu,\Theta))$ 
 on $X^\wedge$, which is introduced in Section 8.2.5 of \cite{EgSc}), the 
 spectrum of $\widehat{A}_{p,q}$ is independent of $1<q<\infty$. Thus we 
 can set $q=2$ and write $\widehat{A}_p=\widehat{A}_{p,2}$. We can 
 assume $p \ge 2$, by passing to the adjoint (i.e., 
 $\widehat{A}_p^*=\widehat{A}_{p'}$ and $-\gamma_p=\gamma_{p'}$). 
 Then ${\rm ker}(\lambda-\widehat{A}_p)\subset
 {\rm ker}(\lambda-\widehat{A}_2) = \{ 0 \}$, since
 $\calK^{2,\gamma_p+2}_2(X^\wedge) \subset \calK^{2,2}_2(X^\wedge)$ in 
 view of $\gamma_{p} \ge \gamma_{2} = 0$. The fact that 
 $\sigma_{M}^{2}(A)(z)$ is invertible for $0 < \re z \le 
 \frac{n+1}{2} - \gamma_{p'} - 2$ implies that
 \begin{equation}\label{kernel}
     {\rm ker}(\lambda-\widehat{A}_{p'}) \subset 
     \calK^{2,2}_2(X^\wedge) = \calK^{2,\gamma_{2}+2}_2(X^\wedge),
     \hspace{0.3cm} \lambda \notin \rpbar;
 \end{equation}
 we shall give the argument below. As a consequence, we have for the 
 adjoint
 $$
 {\rm ker}(\lambda-\widehat{A}_{p})^{*} = 
 {\rm ker}(\bar{\lambda}-\widehat{A}_{p'})  
 \subset 
 {\rm ker}(\bar{\lambda}-\widehat{A}_{2}), 
 \hspace{0.3cm} \lambda \notin \rpbar,
 $$
 hence $\lambda-\widehat{A}_{p}$ is bijective for $\lambda \notin \rpbar$.
 
 In order to see \eqref{kernel} set $\gamma^{1} = 
 {\mathrm{min}}(\gamma_{p'}+2,0)$. The invertibility of the conormal 
 symbol implies that $\lambda-\widehat{A}$ is elliptic with respect 
 to $\gamma^{1}+2$. Moreover, the minimal and maximal extensions of 
 $\lambda-\widehat{A}$, considered as unbounded operators in 
 $\calK^{2,\gamma^{1}}_2(X^\wedge)$, coincide and 
 their domain is $\calK^{2,\gamma^{1}+2}_2(X^\wedge)$. 
 In particular, ${\mathcal N} = 
 {\rm ker} \{ \lambda-\widehat{A} : \calK^{2,\gamma_{p'}+2}_2(X^\wedge)
 \rightarrow \calK^{2,\gamma_{p'}}_2(X^\wedge) \}$ is a subset 
 of the maximal domain, thus it is included in 
 $\calK^{2,\gamma^{1}+2}_2(X^\wedge)$. Iterating this process, we see 
 that ${\mathcal N} \subset \calK^{2,\gamma^j+2}_2(X^\wedge)$ for all 
 $\gamma^j := {\rm min}(\gamma^{j-1}+2,0) = {\rm min}(\gamma_{p'}+2j,0)$. 
 Choosing $j$ large enough we get \eqref{kernel}.
\end{proof}
As a consequence of Theorem \ref{example} we get the following result on the 
maximal regularity for solutions of the Cauchy problem for the Laplacian:  
\begin{theorem}\label{cauchy}
 Let $\Delta$ be the Laplacian on $\intb$ as described above, $1<p<\infty$, 
 and $2\max(p,p')<\dim\bz$. Then the Cauchy problem 
  \begin{equation}\label{cauchy1}
   \dot u(\tau)-\Delta u(\tau)=f(\tau),\quad0<\tau <T;\qquad u(0)=0,
  \end{equation}
 has  a unique solution 
  $$u\in W^1_r\left([0,T],\calH^{0,\gamma_p}_q(\bz)\right)\,\cap\,
    L_r\left([0,T],\calH^{2,\gamma_p+2}_q(\bz)\right)$$
for every
  $$f\in L_r\left([0,T],\calH^{0,\gamma_p}_q(\bz)\right),\qquad 
    1<q,r<\infty.$$ 
 Furthermore, $u$, $u'$, and $\Delta u$ depend continuously on $f$. 
\end{theorem} 
In fact, in Theorem \ref{example} above, we have shown that $-\Delta$ is 
elliptic with respect to $\gamma_p+2$ and any sector $\Lambda_\Delta$ not 
containing $\rpbar$. Moreover, the problem \eqref{cauchy1} is equivalent to 
$\dot v(\tau)-(\Delta-c)v(\tau)=e^{c\tau}f(\tau)$,
$v(0)=0$, and, for sufficiently large 
$c$, the operator $-\Delta + c$ satisfies the assumptions of Theorem 
\ref{bip} for any fixed $0<\theta<\frac{\pi}{2}$ and $\delta > 0$. Then 
 $$\|(-\Delta+c)^{iy}\|_{\calL(\calH^{0,\gamma_p}_q(\bz))}\le 
   c_{p,q}\,e^{\theta|y|}\qquad\forall\;y\in\rz,$$
and Theorem \ref{cauchy}
immediately follows from Theorem 3.2 of \cite{DoVe}. 
\begin{remark} 

{\rm a)}\ 
An interpolation result of Amann \cite{Amann}, Theorem {\rm III.4.10.2},
shows that
 $$W^1_r([0,T],\calH^{0,\gamma_p}_q(\bz))\cap 
   L_r([0,T],\calH^{2,\gamma_p+2}_q(\bz))\hookrightarrow 
   \calC([0,T], (\calH^{0,\gamma_p}_q(\bz),
\calH^{2,\gamma_p+2}_q(\bz))_{\frac{1}{r},r}).$$
We will be interested in the special case, where $r=q$. 
Here, we know from \cite{CSS2}, Corollary {\rm 5.5}, that 
 \begin{equation}\label{interpolation}
   (\calH^{0,\gamma_p}_q(\bz),\calH^{2,\gamma_p+2}_q(\bz))_{\frac{1}{q},q}
   \hookrightarrow \calH^{s,\delta}_q(\bz)\qquad
   \text{for any }
   \mbox{$s<\frac{2}{q^\prime},\,\delta<\gamma_p+\frac{2}{q^\prime}$}.
 \end{equation} 
In particular, we deduce that the solution $u$ in Theorem {\rm \ref{cauchy}}
is a continuous function
on $[0,T]$ with values in $\calH^{s,\delta}_q(\bz)$.

{\rm b)}\ It follows from \cite{Amann}, Theorem {\rm III.4.10.7} and Remark
{\rm III.4.10.9(c)}, that Theorem {\rm\ref{cauchy}}
also holds for initial values 
$u(0) = u_0\in 
(\calH^{0,\gamma_p}_q(\bz),
\calH^{2,\gamma_p+2}_q(\bz))_{\frac{1}{r},r}$.
\end{remark}

Condition \eqref{dimension} implies that $\dim \bz>4$. 
Its strength lies in the fact that it ensures that $\Delta$ is essentially
self-adjoint on $L_2(\bz)$ with domain $\calH^{2,2}_2(\bz)$. 
It was shown by Cheeger that essential selfadjointness also
 holds for $\dim\bz=4$. In that case,
however, the domain is larger than $\calH^{2,2}_2(\bz)$.  For $\dim\bz<4$, the
Laplacian is {\em not} essentially self-adjoint. We then have $\Delta_{\rm min}
\subset \Delta_F\subset\Delta_{\rm max}$, where $\Delta_F$ is Friedrich's
extension. Hence the resolvent set on $L_2(\bz)$ of both $\Delta_{\rm min}$ and
$\Delta_{\rm max}$ is empty, and Theorem \ref{cauchy} will certainly not
be true for the minimal or the maximal extension.


\section{Application: A quasilinear diffusion equation}\label{sec:quasi}

As explained for example in the introduction of \cite{Amann}, diffusion
processes are governed by a quasilinear equation of the form 
$$\dot u(\tau) -\dvz D(u)\grad u(\tau) =f(\tau, u), \quad 0<\tau<T.$$
We now want to illustrate how the above analysis of the Laplacian allows us to solve problems of
this kind for certain choices of the `diffusion matrix' $D$  and the nonlinearity
$f$.

To this end we shall make use of results
obtained in Section 5 of \cite{CSS2} and  a theorem of Cl\'ement and Li, which reads as follows.

\begin{theorem}\label{ClementLi} {\em (Cl\'ement\&Li)}
Let $E_0$ and $E_1$ be Banach spaces, $E_1$ densely and continuously embedded in  $E_0$. For fixed $1<q<\infty$
denote by $E=(E_1,E_0)_{{1}/{q},q}$ the real interpolation space. 

For the quasilinear equation 
 $$\partial_\tau u(\tau) -\tilde{A}(u)u(\tau)=\tilde{f}(\tau,u)+g(\tau)
\ \ on \ ]0,T],\qquad u(0)=u_0\in E,$$
to have a unique solution 
$u\in W^1_q([0,T_1],E_0)\,\cap\,L_q([0,T_1],E_1)$ 
for some $0<T_1\le T$, it is sufficient that there exists an open neighborhood 
$U\subset E$ of $u_0$ such that 
 \begin{itemize}
  \item[(H1)] $\tilde{A}:U\to\calL(E_1,E_0)$ is Lipschitz continuous and $\tilde{A}(u_0)$ 
   is of maximal regularity with respect to $E_0$, $E_1$, and $q$, 
  \item[(H2)] $\tilde{f}:{[0,T]}\times U\to E_0$ is Lipschitz continuous,
  \item[(H3)] $g\in L_q([0,T],E_0)$.
\end{itemize}
\end{theorem}

In (H1), maximal regularity means that the Cauchy problem $\partial_\tau
v-\tilde{A}(u_0)v=h$, $v(0)=v_0$,
has a unique solution $v\in W^1_q([0,T],E_0)\,\cap\,L_q([0,T],E_1)$ 
for every
$h\in L_q([0,T],E_0)$, $v_0\in E$, with $v$, $\partial_t v$, and
$\tilde{A}(u_0)v$ depending continuously on $h$ and $v_0$. 

In order to apply this theorem to our situation, we fix a boundary defining function $t$, which we also use as a coordinate in a neighborhood of the boundary, and choose a Riemannian metric $h_{cone}$ on $\bz$ with a conical degeneracy at the boundary.
We let $\dvz$ and $\grad$ denote the divergence and gradient, respectively, with
respect to $h_{cone}$. More explicitly, writing $h_{cone} = dt^2 +t^2g(t)$ in a neighborhood of $\partial \bz$ with a smooth family $g(\cdot)$ of metrics on the cross-section, we have
$$\grad u = t^{-2}\Big(t^2\partial_tu\
\partial_t +  \sum_{i,j=1}^n g^{ij}\partial_{x_j}u\ \partial_{x_i}\Big).$$
  
Next, we let $a\in\ci(\rz^2)$ denote an arbitrary smooth, positive function. 
We shall consider the case where the diffusion matrix is a scalar multiple of the identity on $T\bz$ 
of the form $D(u)= a(t^{c} u)I_{T\bz}$ with the above boundary defining function $t$ and  an arbitrary positive constant $c$. Here, we identify the complex values of $u$ with elements of $\rz^2$. Instead of being constant, $c$ might be a smooth, real-valued function on $\bz$ which is positive and constant at the boundary. 

We  let 
$E_0=\calH^{0,\gamma_p}_q(\bz)$ and $E_1=\calH^{2,\gamma_p+2}_q(\bz)$ with $p$
and $q$ to be determined later on.
Moreover, we define the second order differential operator
 \begin{equation}\label{divgrad} 
  A(u)=\text{div}(a(t^{c}u)\,\text{grad})
 \end{equation} 

Note that the point evaluation is defined for $u\in
E=(E_1,E_0)_{1/q,q}$ by the Sobolev embedding theorem provided $q>(n+3)/2$. 
We are going to show the following theorem: 
\begin{theorem}\label{quasi}
Assume that $\dim \bz>4$ and let $T>0$. 
Then there exists a suitable choice of numbers $1<p,q<\infty$ and  
$0<T_1\le T$ such that the equation  
\begin{equation}\label{quasieq} 
\partial_\tau u(\tau) -A(u)u(\tau)= f(\tau,u)+g(\tau),\qquad u(0)=u_0\in E, 
\end{equation} 
has a unique solution 
  $$u\in W^1_q([0,T_1],E_0)\,\cap\,L_q([0,T_1],E_1)$$
 for every initial value $u_0\in\cicomp(\intb)$, 
every $f\in {\rm Lip}([0,T]\times U,E_0)$,
and every $g\in L_q([0,T],E_0)$.
\end{theorem}

As we showed in \cite{CSS2}, Corollary {\rm 5.11}, examples of  functions $f$
satisfying the above assumption $($for suitable $p$, $q)$ include $f(u) = |u|^\alpha$, $\alpha \in\rz$, or $f(u) = u^\alpha$, $\alpha\ge 1$ $($and hence their linear combinations$)$. A specific example
here is the time-dependent Ginzburg-Landau equation
$\partial_\tau u -\Delta u = u-u^3 $
for initial data $u(0) =u_0\in \cicomp(\intb)$. 

We deduce Theorem \ref{quasi} from \ref{ClementLi}.   
As a preparation we rewrite equation \eqref{quasieq} as 
\begin{equation}\label{umf}
\partial_\tau u -\tilde{A}(u)u=\tilde{f}(u)+g,\quad u(0)=u_0\in E
\end{equation}
with $\tilde{A}(u)=a(t^{c}u)\,\Delta$ and
\begin{equation}
\label{ftilde}
  \tilde{f}(\tau,u)=f(\tau,u)-
  (\partial_1a)(t^cu)\skp{\grad(t^c\Re u)}{\grad u}-
i (\partial_2a)(t^cu)\skp{\grad(t^c\Im u)}{\grad u}.
\end{equation}
Here, $\skp{\cdot}{\cdot}$ is the complexified pointwise scalar product on
$T\bz$ given by $h_{cone}$ and $\Delta$ is the associated Laplacian. 
For large $q$, solving \eqref{umf} in 
$W^1_q([0,T_1],E_0)\,\cap\,L_q([0,T_1],E_1)$
is equivalent to solving \eqref{quasieq} in that space. 
In fact, for $q>n+3$ the Sobolev embeddng theorem implies that $E\subseteq
C^1(\intb)$. Hence the solutions to both equations will be functions continuous
in $\tau$ and continuously differentiable in $(t,x)$ so that $A(u)u+f(\tau,u)$
and $\tilde A(u)u + \tilde f(\tau,u)$ coincide as distributions on $\intb$ for
each $\tau$.

In \cite{CSS2}, Theorem 5.7, we already have shown that $\tilde{A}$ satisfies
condition (H1), provided $p<\frac{n+1}{2}$ is close to
$\frac{n+1}{2}$, and $q$ is large.  
It is therefore sufficient to show the Lipschitz continuity of the map 
 \begin{equation}\label{abcde}
  u\mapsto (\partial_1a)(t^cu)\skp{\grad(t^c\Re u)}{\grad u}:U\to E_0
 \end{equation}
with $p$ and $q$ subject to the above condition; the other term  in \eqref{ftilde} can be treated in the same way. 

\begin{lemma}Let $c>0$. For $p<\frac{n+1}{2}$ sufficiently close to $\frac{n+1}{2}$ and $q$ sufficiently large, the mapping 
$u\mapsto (\partial_1a)(t^cu):U\to L_\infty(\bz)$ is Lipschitz continuous on bounded subsets of $E$.
\end{lemma}

\begin{proof}In \cite{CSS2}, Lemma 5.6, it is shown that then $E\hookrightarrow
t^c\calC(\bz)$,
where the right-hand side denotes the space of bounded continuous functions on
$\bz$, multiplied by $t^c$.
Thus
$$\|(\partial_1a)(t^cu_1)- (\partial_1a)(t^cu_2)\|_{L_\infty(\bz)}
\le C\|t^cu_1-t^cu_2\|_{L_\infty(\bz)}\le C\|u_1-u_2\|_E,$$ 
where the first constant $C$ is the maximum of the norm of the total derivative of $\partial_1a(s)$, as $s$ varies over the bounded set of all values of $t^cu$, $u\in U$. \end{proof}

In order to finish the proof of Theorem
\ref{quasi} it then is enough to establish the Lipschitz continuity of
$$u\mapsto \skp{\grad(t^c\Re u)}{\grad u}= t^c\skp{\grad(\Re u)}{\grad u}+ \skp{\Re u\grad(t^c)}{\grad u}.$$
We observe that the right-hand side is a linear combination of terms of the
form $(D_1u)(D_2u)$, where 
$$D_j: \calH^{s,\delta}_q(\bz)\to\calH^{s-1,\delta-1+c/2}_q(\bz)$$
is a  bounded real-linear operator for all $s$, $\delta$, and $q$. In fact, 
close to the boundary we have the identities 
$$t^c \skp{\grad \Re u}{\grad u}= (t^{c/2-1} t\partial_t\Re u)(t^{c/2-1}
t\partial_t u)+\sum_{i,j} g^{ij}
(t^{c/2-1}\partial_{x_i}\Re u)(t^{c/2-1}\partial_{x_j} u)$$
and
$$ \Re u\skp{\grad t^c}{\grad u} = c\ \Re u \ t^{c-1}\partial_t
u = c(t^{c/2-1}\Re
u)(t^{c/2-1} t\partial_t u).
$$

The following rather technical proposition, in connection with Corollary \ref{lipschitz}, below, treats the Lipschitz continuity for these functions:
\begin{proposition}
 Let $D_1$, $D_2$ be as above, $1<q<\infty$, and $s,\delta,\gamma\in\rz$ be 
 such that $s-1>\frac{n+1}{q}$ and $\delta-1+c/2>(n+1+2\gamma)/4$.
Then the map  
  $$u\mapsto (D_1u)(D_2u):
    \calH^{s,\delta}_q(\bz)\longrightarrow \calH^{0,\gamma}_q(\bz)$$ 
 is Lipschitz continuous on bounded sets. 
\end{proposition}
\begin{proof}
The proof is similar to that of \cite{CSS2}, Theorem 5.15.
For the convenience of the reader we provide the details. Choose a cut-off
function $\omega\in \cicomp([0,1[)$ and $\psi\in \cicomp(\intb)$ with
$\omega^2+\psi^2=1$. Then $(D_1u)(D_2u)=(\omega D_1u)(\omega D_2u)+(\psi
D_1u)(\psi D_2u)$. 
We first focus on the analysis near the boundary, i.e.\ we show that 
$$u \mapsto (\omega D_1u)(\omega D_2u):
\calH^{s,\delta}_q(\bz )\to\calH^{0,\gamma}_q(\bz)$$
is Lipschitz continuous on bounded sets. 
%
We abbreviate $u_j=D_ju$ and $v_j=D_jv$, $j=1,2$, for $u,v\in
\calH^{s,\delta}_q(\bz)$. Since $u_j$ and $v_j$ have their support near
$t=0$, we have according to \eqref{9.5}  
\begin{align*}
\|u_1u_2&-v_1v_2\|_{\calH^{0,\gamma}_q(\bz)}=
\|S_\gamma(u_1u_2)-S_\gamma(v_1v_2)\|_{L_q(\rz\times X)}
=\|(S_{\tilde\gamma}u_1) (S_{\tilde\gamma}u_2) -( S_{\tilde\gamma}v_1)(S_{\tilde\gamma}v_2)\|_{L_{q}(\rz\times
X)}\\
&\le\| S_{\tilde\gamma}u_1\|_{L_{2q}(\rz\times X)}\|S_{\tilde\gamma}(u_2-v_2)\|_{L_{2q}(\rz\times X)} + 
\| S_{\tilde\gamma}v_2\|_{L_{2q}(\rz\times X)}\|S_{\tilde\gamma}(u_1-v_1)\|_{L_{2q}(\rz\times X)}. 
\end{align*}
Here, we chose $\tilde \gamma = (n+1+2\gamma)/4$, and we employed H\"older's inequality. Next we use the embedding 
$H^{s-1}_q(\rz^{n+1})\hookrightarrow L_{2q}(\rz^{n+1})$, valid for 
$s-1-\frac{n+1}q \ge -\frac{n+1}{2q}$, cf.\ \cite{Triebel}, 2.8.1 Remark 2. 
Since we assumed that $s-1>\frac{n+1}q$
we deduce that 
$$
{\| S_{\tilde\gamma}u_1\|_{L_{2q}(\rz\times X)}\le \|
S_{\tilde\gamma}u_1\|_{H^{s-1}_{q}(\rz\times X)}
\le C\ \|u_1\|_{\calH^{s-1,\tilde\gamma}_q(X^\wedge)}}\le  C \
\|u\|_{\calH^{s,\tilde\gamma+1-c/2}(\bz)}\le
 C \ \|u\|_{\calH^{s,\delta}_q(\bz).}
$$
The second estimate results from the continuity of $D_1$; for the third we used that $\tilde \gamma \le
\delta-1+c/2$.
In the same way we estimate  $ \| S_{\tilde\gamma}v_2\|_{L_{2q}(\rz\times X)}$  and finally 
$$\|S_{\tilde\gamma}(u_j-v_j)\|_{L_{2q}(\rz\times X)}\le  C
\|u-v\|_{\calH^{s,\delta}_q(\bz)},\ \ j=1,2.$$
Next set  $u_j=\psi D_ju$ and $v_j=\psi
D_ju$ and note that the norm of
$\calH^{0,\gamma}_q(\bz)$ coincides with that of $L_q(\bz)$ on their
supports. Then the estimate 
$$\|u_1u_2-v_1v_2\|_{L_q(\bz)}\le
\|u_1\|_{L_{2q}(\bz)}\|u_2-v_2\|_{L_{2q}(\bz)}
+\|v_2\|_{L_{2q}(\bz)}\|u_1-v_1\|_{L_{2q}(\bz)}$$
plus the fact that the norms of $u_1$, $v_2$ and $u_j-v_j$ can be estimated  by
the norms of $u$, $v$, and $u-v$ in $\calH^{\delta,\gamma}_q(\bz)$ completes
the argument.
\end{proof} 


\begin{corollary}\label{lipschitz} 
 Let $D_1$, $D_2$ as above. Then there exist $1<p,q<\infty$ such that the map  
  $$u\mapsto (D_1u)(D_2u):E\longrightarrow E_0=\calH^{0,\gamma_p}_q(\bz)$$
 is Lipschitz continuous on bounded sets.   
\end{corollary}
\begin{proof}
 We have $E\hookrightarrow\calH^{s,\delta}_q(\bz)$ for any 
 $s<{2}/{q^\prime}$ and $\delta<\gamma_p+{2}/{q^\prime}$. 
 Choosing  $p<({n+1})/2$ close to  $({n+1})/2$ and $q$ large, we have 
$$ \gamma_p+\frac2{q'}-1+\frac c2>
\frac{n+1+2\gamma_p}4\ \ \ {\rm and } \ \ \ \frac2{q'}-1>\frac{n+1}q .$$
We can then pick $s,\delta$ with $ \gamma_p+\frac2{q'}-1+\frac
c2>\delta-1+\frac c2>
\frac{n+1+2\gamma_p}4$ and $
\frac2{q'}-1>s-1>\frac{n+1}q $ 
and  apply Lemma \ref{lipschitz}. 
\end{proof}

\section{Appendix: Smoothing Mellin symbols and Green symbols}
\label{appendix}
The structure of the resolvent (respectively parametrix) of a differential 
operator $A$ as given in Theorem \ref{parametrix} at the first glance does 
not coincide with those which you find for example in \cite{EgSc} or 
\cite{Gil}. This is mainly due to the fact that we consider $A$ as an 
unbounded operator in $\hsgpb$ whose resolvent acts continuously in 
$\hsgpb$, and do not consider $A$ as a bounded operator acting from 
$\calH^{s+\mu,\gamma+\mu}_p(\bz)$ to $\hsgpb$. We shall use this appendix 
to clarify this point. 

Let us begin with a discussion of so-called {\em Green symbols}. Let us set 
 $$\calK^{s,\gamma}(X^\wedge)^\nu=
   \spk{t}^\nu\calK^{s,\gamma}_2(X^\wedge)$$ 
for real $\nu$, cf.\ Definition \ref{ksgp}. These are Hilbert spaces, 
and $\calK^{-s,-\gamma}(X^\wedge)^{-\nu}$ can be identified with the 
dual space of $\calK^{s,\gamma}(X^\wedge)^\nu$ via the scalar-product 
in $\calK^{0,0}(X^\wedge)$. The operators $\kappa_\varrho$ defined in 
\eqref{kappa} extend by continuity to operators in 
$\calL(\calK^{s,\gamma}(X^\wedge)^\nu)$. 

For $\mu\in\rz$ and $d>0$ we let 
 $$S^{\mu,d}(\Lambda;\calK^{s,\gamma}(X^\wedge)^\nu,
   \calK^{s',\gamma'}(X^\wedge)^{\nu'})$$
denote the space of all smooth functions 
$a\in\ci(\Lambda,\calL(\calK^{s,\gamma}(X^\wedge)^\nu,
\calK^{s',\gamma'}(X^\wedge)^{\nu'}))$ satisfying
 $$\|\kappa_{\spk{\lambda}^{-1/d}}\{\partial_\lambda^\alpha a(\lambda)\}
   \kappa_{\spk{\lambda}^{1/d}}\|\le 
   c_\alpha\,\spk{\lambda}^{\frac{\mu}{d}-|\alpha|}$$
uniformly for $\lambda\in\Lambda$ and all multiindices $\alpha$. 

We call a smooth function  
$b\in\ci(\Lambda_\Delta\setminus0,\calL(\calK^{s,\gamma}(X^\wedge)^\nu,
\calK^{s',\gamma'}(X^\wedge)^{\nu'}))$ 
{\em twisted homogeneous} of degree $(\mu,d)$ if it fulfills  
 $$b(\varrho^d\lambda)=\varrho^\mu\,\kappa_\varrho\,b(\lambda)\,
   \kappa_\varrho^{-1}$$
for all $\lambda$ and $\varrho>0$. Note that multiplying $b$ with a 
0-excision function (supported sufficiently far away from zero) yields a 
symbol in $S^{\mu,d}(\Lambda;\calK^{s,\gamma}(X^\wedge)^\nu,
\calK^{s',\gamma'}(X^\wedge)^{\nu'})$. The space 
 $$S^{\mu,d}_{cl}(\Lambda;\calK^{s,\gamma}(X^\wedge)^\nu,
 \calK^{s',\gamma'}(X^\wedge)^{\nu'})$$
then consists of all symbols from 
$S^{\mu,d}(\Lambda;\calK^{s,\gamma}(X^\wedge)^\nu,
\calK^{s',\gamma'}(X^\wedge)^{\nu'})$, that have asymptotic expansions 
$a\sim\sum\limits_{k=0}^\infty a^{(\mu-j,d)}$ with functions 
$a^{(\mu-j,d)}$ that are twisted homogeneous of degree $(\mu-j,d)$. 
\begin{definition}\label{green3}
 Let $\gamma,\gamma'\in\rz$. If $g\in 
 S^{\mu,d}(\Lambda;\calK^{0,\gamma}(X^\wedge),\calK^{0,\gamma'}(X^\wedge))$, 
 we can form the adjoint symbol $g^*\in
 S^{\mu,d}(\Lambda;\calK^{0,-\gamma'}(X^\wedge),\calK^{0,-\gamma}(X^\wedge))$ 
 by taking pointwise the adjoint with respect to the 
 $\calK^{0,0}(X^\wedge)$-scalar product. We then call $g$ a Green symbol if 
 additionally there exists an $\eps=\eps(g)>0$, such that  
  \begin{align*}
   g&\in\mathop{\mbox{\Large$\cap$}}_{s,s',\nu,\nu'} 
   S^{\mu,d}_{cl}(\Lambda;\calK^{s,\gamma}(X^\wedge)^\nu,
   \calK^{s',\gamma'+\eps}(X^\wedge)^{\nu'}), \\
   g^*&\in\mathop{\mbox{\Large$\cap$}}_{s,s',\nu,\nu'} 
   S^{\mu,d}_{cl}(\Lambda;\calK^{s,-\gamma'}(X^\wedge)^\nu,
   \calK^{s',-\gamma+\eps}(X^\wedge)^{\nu'}).
  \end{align*}
 The entity of all such Green symbols we shall denote by 
  $$R^{\mu,d}_G(X^\wedge;\Lambda,(\gamma,\gamma')).$$
\end{definition}
It is a trivial fact that if $\gamma'\ge\gamma''$, then 
  $$R^{\mu,d}_G(X^\wedge;\Lambda,(\gamma,\gamma'))\subset
    R^{\mu,d}_G(X^\wedge;\Lambda,(\gamma,\gamma'')).$$
Moreover it can be shown, cf.\ \cite{Seil}, that in case 
$\gamma=\gamma'$ both Definitions \ref{green3} and \ref{green2} yield the 
same symbols respectively operator-families. In other words, Green symbols 
can either be characterized by their mapping properties in Sobolev spaces 
or by the structure of their kernels. 

Let us now return to the resolvent, cf.\ Theorem \ref{parametrix}. If you 
compare with \cite{Gil}, you will find that there `our' term $G(\lambda)$ 
is replaced by a term of the form $G_0(\lambda)+M(\lambda)$, where 
 $$G_0(\lambda)\in
   R^{-\mu,\mu}_G(X^\wedge;\Lambda,(\gamma,\gamma+\mu)),$$
is a Green symbol and 
 $$M(\lambda)=\omega(t[\lambda]^{\frac{1}{\mu}})\,t^\mu\,
   \opm{\gamma-\frac{n}{2}}(h)\,\omega_0(t[\lambda]^{\frac{1}{\mu}})$$
for some cut-off functions $\omega,\omega_0\in\cicomp(\rpbar)$, and a 
meromorphic Mellin symbol 
 $$h\in M^{-\infty}_P(X).$$
The last notation roughly means that $h$ is a meromorphic function on the 
complex plane with values in $L^{-\infty}(X)$,
the smoothing pseudodifferential operators on 
$X$, having only finitely many poles in any vertical 
strip $|\re z|\le\beta$, $\beta>0$, and the Laurent coefficients of the 
principal part of $h$ at any such pole are finite rank operators. For more 
details see \cite{EgSc}, Section 8.1.2. By the above observation, 
 $$G_0(\lambda)\in
   R^{-\mu,\mu}_G(X^\wedge;\Lambda,\gamma).$$
The same is also true for $M$, since it is easy to see that 
 $$M\in 
   S^{-\mu,\mu}_{cl}(\Lambda;\calK^{s,\gamma}(X^\wedge)^\nu,
   \calK^{s',\gamma+\mu}(X^\wedge)^{\nu'})$$
for all $s,s',\nu,\nu'$ (note that $M$ is twisted homogeneous for large 
$|\lambda|$). Observe that $\mu \ge 1$ since we are dealing with 
differential operators of positive order. The adjoint symbol is given by 
 $$M^*(\lambda)=\omega_0(t[\lambda]^{\frac{1}{\mu}})\,t^\mu\,
   \opm{-\gamma_+-\frac{n}{2}}(h^*)\,\omega(t[\lambda]^{\frac{1}{\mu}})\;+\;
   G_1(\lambda)$$
where $h^*(z)=h(n+1-\mu-\overline{z})^*$ and  
 $$G_1(\lambda)=t^\mu\,\omega_0(t[\lambda]^{\frac{1}{\mu}})\,
   \Big\{\opm{-\gamma-\mu-\frac{n}{2}}(h^*)-
   \opm{-\gamma_+-\frac{n}{2}}(h^*)\Big\}\,
   \omega(t[\lambda]^{\frac{1}{\mu}}).$$
Here, $\gamma_+=\gamma$ if $h^*$ has no pole on the line 
$\re z=\frac{n+1}{2}+\gamma$, otherwise $\gamma_+>\gamma$ sufficiently 
close to $\gamma$. However, it is known, cf.\ \cite{EgSc}, Section 
8.1.2, Theorem 6, that then 
$G_1\in R^{-\mu,\mu}_G(X^\wedge;\Lambda,(-\gamma, -\gamma))$ and 
 $$M^*-G_1\in
   S^{-\mu,\mu}_{cl}(\Lambda;\calK^{s,-\gamma}(X^\wedge)^\nu,
   \calK^{s',-\gamma_++\mu}(X^\wedge)^{\nu'})$$
for all $s,s',\nu,\nu'$. All together this shows that 
$M\in R^{-\mu,\mu}_G(X^\wedge;\Lambda,\gamma)$ and hence justifies the 
description of the resolvent we have given in Theorem \ref{parametrix}. 

\section{Notation}\label{notation} 
For $0\not=\lambda\in\cz$ we let $\arg\lambda$ be the unique number 
$-\pi\le\arg\lambda<\pi$ such that $\lambda=|\lambda|e^{i\arg\lambda}$. For 
$z\in\cz$ we then set 
 $$
  \lambda^z=|\lambda|^ze^{iz\arg\lambda}.
 $$
For fixed $z$ this is a holomorphic function in 
$\lambda\in\cz\setminus\{\lambda\in\rz\st\lambda\le0\}$. 

For $\delta>0$ and $0<\theta<\pi$ we let $\Lambda=\Lambda(\delta,\theta)$ 
denote the closed {\em keyhole region} 
 $$
  \Lambda(\delta,\theta)=\{\lambda\in\cz\st|\lambda|\le\delta
                           \text{ or }|\arg\lambda|\ge\theta\}
 $$
and $\calC=\calC(\delta,\theta)$ its parametrized boundary, 
$\calC=\calC_1\cup\calC_2\cup\calC_3$, with 
 \begin{equation}\label{not3}
  \calC_1(t)=te^{i\theta},\;-\infty<t\le\delta;\quad
  \calC_2(t)=\delta e^{-i t},\;-\theta\le t\le\theta;\quad
  \calC_3(t)=te^{-i\theta},\;\delta\le t<\infty. 
 \end{equation}
We let $\Lambda_\Delta=\Lambda_\Delta(\theta)$ denote the closed sector 
contained in $\Lambda$, 
 \begin{equation}\label{not4}
  \Lambda_\Delta(\theta)=\{\lambda\in\cz\st|\arg\lambda|\ge\theta\}\cup\{0\}
 \end{equation}
and, similar to \eqref{not3}, $\calC_\Delta$ its parametrized boundary. 

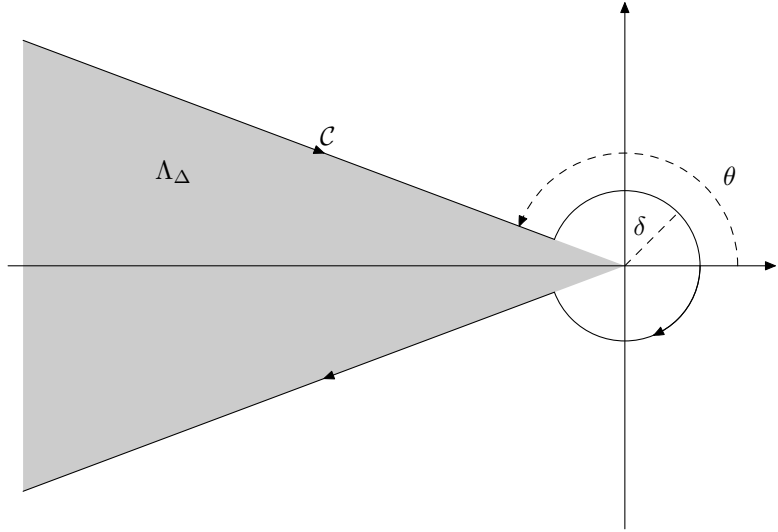
\begin{center}
    \includegraphics{figures.1}
\end{center}

We now recall various spaces of pseudodifferential symbols and operators we 
shall use throughout this paper. In the following we let $\mu,d\in\rz$ and 
$d$ positive. 

We call a function smooth on $\Lambda$, if it is the restriction to 
$\Lambda$ of a function which is smooth in an open neighborhood of 
$\Lambda$. If $E$ is a Fr\'echet space, then 
 \begin{equation}\label{not5}
  \calS(\Lambda,E)
 \end{equation}
consists of all $u\in\ci(\Lambda,E)$ satisfying 
 $$\sup_{\lambda\in\Lambda}\trinorm{\partial^\gamma_\lambda u(\lambda)}\,
   |\lambda|^N<\infty$$
for any multi-index $\gamma\in\nz_0^2$, any $N\in\nz$, and any continuous 
semi-norm $\trinorm{\cdot}$ of $E$. The space of symbols of order $\mu$ and 
anisotropy $d$, 
 $$
   S^{\mu,d}(\rz^m_y\times\rz^n_\eta;\Lambda),
 $$
consists of all functions (possibly matrix-valued)
$a\in\ci(\rz^m\times\rz^n\times\Lambda)$, which 
fulfill the estimates 
 $$|\partial_y^\beta\partial_\eta^\alpha\partial_\lambda^\gamma 
   a(y,\eta,\lambda)|\le c_{\alpha\beta\gamma}\,
   \spk{\eta,\lambda}^{\mu-|\alpha|-d|\gamma|}_d,\qquad
   \spk{\eta,\lambda}_d=(1+|\eta|^2+|\lambda|^{\frac{2}{d}})^{\frac{1}{2}},$$
for any multi-indices $\alpha$, $\beta$, and $\gamma$. Further we set 
 $$
  S^{\mu,d}(\rpbar\times\rz^{m-1}\times\rz^n;\Lambda)=
  S^{\mu,d}(\rz^m\times\rz^n;\Lambda)|_{\rpbar\times\rz^{m-1}\times
  \rz^n\times\Lambda}. 
 $$
For a compact manifold $X$, $\text{\rm dim}\,X=n$, the space 
 \begin{equation}\label{not6}
  L^{\mu,d}(X;\Lambda) 
 \end{equation}
of parameter-dependent pseudodifferential operators of order $\mu$ and 
anisotropy $d$ (acting on sections of a vector bundle) consists
of all operator-families, which are obtained as a 
sum (according to a covering of $X$ by coordinate neighborhoods) of local 
operators with symbols from $S^{\mu,d}(\rz^n\times\rz^n;\Lambda)$ 
and a smoothing remainder from 
$L^{-\infty}(X;\Lambda):=\calS(\Lambda,L^{-\infty}(X))$. In the last 
definition, $L^{-\infty}(X)$ is the usual space of smoothing operators on 
$X$, i.e.\ the space of all integral operators having a smooth kernel. 

If $\gamma \in \rz$ and $\Gamma_{\gamma}$ denotes the vertical line in 
the complex plane
$$
\Gamma_{\gamma} = \{ z \in \cz \; | \; \re z = \gamma \},
$$
the space of symbols 
 $$
  MS^{\mu}(\rz_+\times\rz^n\times
   \Gamma_{\frac{n+1}{2}-\gamma}\times\rz^n) 
 $$ 
consists of all functions $a\in\ci(\rz_+\times\rz^n\times    
\Gamma_{\frac{n+1}{2}-\gamma}\times\rz^n)$ which satisfy the estimates
 $$|\partial^l_{\tau} (t\partial_{t})^k 
     \partial^\alpha_{\xi} \partial^\beta_{x}
     a(t,x,\mbox{$\frac{n+1}{2}-\gamma$}+i\tau,\xi) | 
   \le c_{kl\alpha\beta} \langle\tau,\xi \rangle^{\mu-l-|\alpha|},
      \qquad\langle\tau,\xi \rangle = (1 + \tau^{2} + |\xi|^{2})^{1/2}.$$
The associated (Fourier-Mellin) pseudodifferential operator is 
 $$[\opm{\gamma-\frac{n}{2}}(a)u](t,x)=
   \int_{\re z=\frac{n+1}{2}-\gamma}t^{-z}\op(a)(t,z)(\calM u)(z,x)\,
   \dbar z,
   \qquad u \in \cicomp(\rz_{+}\times \rz^n),$$ 
where $\op$ is the standard Fourier pseudodifferential operator on $\rz^n$, 
and $\calM$ the Mellin transform 
 $$(\calM v)(z)=\int_0^\infty t^z v(t)\,\mbox{$\frac{dt}{t}$}.$$ 
The  operator $\opm{\gamma-\frac{n}{2}}(a)$ induces continuous mappings on the associated Mellin Sobolev spaces 
$$\opm{\gamma-\frac{n}{2}}(a):  
\calH^{s,\gamma}_q(\rz_+\times\rz^n)\to\calH^{s-\mu,\gamma}_q(\rz_+\times\rz^n),$$ 
$s\in\rz$, $1<q<\infty$, where 
 \begin{equation}\label{9.5}
  \calH^{s,\gamma}_q(\rz_+\times\rz^n):=\big\{u\st 
  (Su)(t,x) :=
  e^{(\frac{n}{2}-\gamma)t}u(e^{-t},x)\in H^s_q(\rz^{1+n}_{(r,x)})\big\}.
 \end{equation}
 The continuity is due to the fact that a Mellin pseudodifferential 
 operator on $\rz_{+} \times \rz^n$ transforms, under conjugation by 
 $S$, to a usual pseudodifferential operator on $\rz^{1+n}$, and then 
 the Calder\'on-Vaillancourt theorem applies.

Using local coordiantes, we obtain the spaces  $  \calH^{s,\gamma}_q(X^\wedge)$
for an arbitrary closed manifold $X$ with the corresponding map $S_\gamma :
\calH^{s,\gamma}_q(X^\wedge)\to H^s_q(\rz\times X).$

The space $ \calH^{s,\gamma}_q(\bz)$ consists of all $u\in H^{s}_{q,{\rm loc}
}({\rm int}\,\bz)$  such that
$\omega u\in \calH^{s,\gamma}_q(X^\wedge)$ for a cut-off function $\omega$. 

\setlength{\parskip}{0pt}

\begin{small}
\bibliographystyle{amsalpha}

\end{small}

\end{document}